\documentclass[10pt]{article}
\usepackage{amsmath,amssymb,bbm}
\usepackage{a4wide}
\usepackage{url}
\usepackage[all]{xy}

\def\carac#1,#2{
\left[
\begin{smallmatrix}
#1 \\ #2
\end{smallmatrix}
\right]
}

\newenvironment{prf}[1]{\trivlist 
\item[\hskip \labelsep{\bf 
#1.\hspace*{.3em}}]}{~\hspace{\fill}~$\square$\endtrivlist} 
\newenvironment{proof}{\begin{prf}{Proof}}{\end{prf}} 

\newtheorem{theorem}{Theorem}[section]
\newtheorem{lemma}[theorem]{Lemma}
\newtheorem{proposition}[theorem]{Proposition}
\newtheorem{corollary}[theorem]{Corollary}

\newtheorem{definition}[theorem]{Definition}

\newcommand{\hsp}{\hspace{5pt}}
\newcommand{\proj}[1]{ \mathrm{Proj} (#1)}
\newcommand{\spec}[1]{ \mathrm{Spec} (#1)}
\newcommand{\cA}{\mathcal{A}}
\newcommand{\cM}{\mathcal{M}}

\renewcommand{\AA}{\mathbb{A}}
\newcommand{\xq}{\mathbb{Q}}
\newcommand{\QQ}{\mathbb{Q}}

\newcommand{\xz}{\mathbb{Z}}
\newcommand{\ZZ}{\mathbb{Z}}
\newcommand{\xf}{\mathbb{F}}
\newcommand{\FF}{\mathbb{F}}
\newcommand{\xg}{\mathbb{G}}
\newcommand{\xp}{\mathbb{P}}
\newcommand{\PP}{\mathbb{P}}

\newcommand{\mat}[2]{\mathbb{M}_{#1}(#2)}

\newcommand{\Aut}{\mathrm{Aut}}

\newcommand{\pol}{\mathcal{L}}

\newcommand{\emm}{\mathcal{M}}
\newcommand{\xone}{\mathbbm{1}}

\newcommand{\ignore}[1]{}

\title{Higher dimensional $3$-adic CM construction}
\author{Robert Carls, David Kohel, David Lubicz
\footnote{This work was supported by the
Australian Research Council grant DP0453134.}
\date{}
}

\begin{document}

\maketitle

{\small 
\begin{center}
\begin{tabular}{l@{\hspace{0.5cm}}l}
Robert Carls, David Kohel & David Lubicz \\
\texttt{\{carls,kohel\}@maths.usyd.edu.au} & 
\texttt{david.lubicz@univ-rennes1.fr}\\
School of Mathematics \& Statistics F07 & CELAR \\
University of Sydney NSW 2006 & BP 7419 35174 Bruz Cedex\\
Australia & France \\
\end{tabular}
\end{center}
}

\bigskip
\begin{abstract}
\noindent
We find equations for the higher dimensional analogue of the modular curve 
$X_0(3)$ using Mumford's algebraic formalism of algebraic theta functions.
As a consequence, we derive a method for the construction of genus $2$ 
hyperelliptic curves over small degree number fields whose Jacobian has 
complex multiplication and good ordinary reduction at the prime $3$.  
We prove the existence of a quasi-quadratic time algorithm for computing 
a canonical lift in characteristic $3$ based on these equations, with a  
detailed description of our method in genus $1$ and $2$.
  \\
  \\
  \textbf{Keywords:} CM-methods, canonical lift, theta functions, modular equations.
\end{abstract}

\section{Introduction}

The theory of complex multiplication yields an efficient method to
produce abelian varieties over a finite field with a prescribed
endomorphism ring.  In the case of elliptic curves, one starts with
$\mathcal{O}$ an order in an imaginary quadratic field of discriminant
$D$. Let $h=h(D)$ be the class number of $\mathcal{O}$.  It is well
known~\cite[Ch.II]{MR1312368} that there exist exactly $h$ isomorphism
classes of elliptic curves with complex multiplication by
$\mathcal{O}$. Let $j_i$ be their $j$-invariants where $i=1, \ldots h$
.  The usual CM-method for elliptic curves consists of computing the
$j_i$ using floating point arithmetic.  One then recovers the Hilbert
class polynomial $$H_D(X)=\prod_{i=1}^h(X-j_i)$$ from its real
approximation, using the fact that it has integer coefficients. It is
usual to assess the complexity of this algorithm with respect to the
size of the output.  As $h$ grows quasi-linearly with respect to $D$,
the complexity parameter is $D$.  
\newline\indent 
In 2002, Couveignes
and Henocq~\cite{MR2041087} introduced the idea of CM construction via
$p$-adic lifting of elliptic curves.  The basis of their idea is to
construct a \emph{CM lift}, i.e.~a lift of a curve over a finite field
to characteristic zero such that the Jacobian of the lifted curve has
complex multiplication.  In the ordinary case the lifting can be done
in a canonical way.  In fact, one lifts a geometric invariant of the
curve using modular equations.  The computation of the canonical lift
of an ordinary elliptic curve has drawn a considerable amount of
attention in the past few years following an idea of
Satoh~\cite{MR2001j:11049,MR1895422,LL03,Gaudry2002,SaSkTa2001,Mestre4}.  Mestre
generalized Satoh's method to higher dimension using theta constants.
His purely $2$-adic method \cite{Mestre5,Ritzenthaler03} is based on a
generalization of Gauss' \emph{arithmetic geometric mean}~(AGM)
formulas.  In this article we present formulas which may be seen as a
$3$-adic analogue of Mestre's generalized AGM equations.  In contrast
to the latter ones, our equations do not contain information about the
action of a lift of relative Frobenius on the cohomology.  In order to
construct the canonical lift, we apply a modified version of the
lifting algorithm of Lercier and Lubicz~\cite{ll05} to our equations.
\newline\indent 
Next we compare our $3$-adic CM method to the $2$-adic CM
method for genus $2$ of Gaudry et al.~\cite{Gaudryandall},
which uses the classical Richelot
correspondence for canonical lifting.
The latter method applies only to those CM fields $K$
in which the prime $2$ splits completely in the quadratic extension
$K/K_0$, where $K_0$ is the real subfield.  For any other CM field
$K$, the reduction of the CM curve at $2$ will be non-ordinary. Thus
there exists no ordinary curve with CM by $K$ to serve as input to the
algorithm.  The method presented in this paper exchanges this
condition at $2$ with the analogous condition at $3$.  Hence the
resulting $3$-adic CM method applies to a large class of CM fields
which are not treatable by the prior $2$-adic CM
method~\cite{Gaudryandall}.  \newline\indent Finally, we describe the
techniques that are used in order to prove the equations introduced in
the present paper.  We prove our equations using the theory of
algebraic theta functions which was developed by Mumford~\cite{mu66}.
In the $3$-adic arithmetic situation we make use of a canonical
coordinate system on the canonical lift whose existence is proven
in~\cite{ca05a}.  Our algorithm is proven by $3$-adic analytic means
and Serre-Tate theory.  \newline\indent This article is structured as
follows. In Section~\ref{main} we prove equations which are satisfied
by the canonical theta null points of canonical lifts of ordinary abelian varieties
over a perfect field of characteristic~$3$. 
In Section~\ref{basic}, for lack of a suitable reference, we prove some
properties of algebraic theta functions which are used in the proof
of the modular equations for the prime $3$.
In Section~\ref{algo} we
describe a method for CM construction via canonical lifting of abelian
surfaces in characteristic $3$.
In Section~\ref{riemequations} we recall classical results about the
moduli of hyperelliptic genus $2$ curves and provide
examples of the CM invariants of abelian surfaces and genus 2 curves.

\section{Modular equations of degree $3$ and level $4$}
\label{main}

In this section we prove equations which have as solutions the theta
null points of the canonical lifts of ordinary abelian varieties over
a perfect field of characteristic $3$. The latter equations form an
essential ingredient of the $3$-adic CM construction which is given
in Section \ref{algo}.
Our proof uses Mumford's formalism of \emph{algebraic
theta functions} \cite{mu66}. The results of Section \ref{galcantheta}
cannot be obtained in a complex analytic setting.
We remark that in \cite{ko98} Y. Kopeliovich proves
higher dimensional theta identities of degree $3$ using complex
analytic methods. Our purely algebraic method yields similar
equations. Our set of equations is 'complete' in the sense
that it defines a higher dimensional analogue of the classical
modular curve $X_0(3)$.

\subsection{Theta null points of $3$-adic canonical lifts}
\label{modeq3}

For the basics about algebraic theta functions and standard notation we refer to \cite{mu66}.
Let $R$ be a complete noetherian local ring with perfect residue field $k$
of characteristic $3$. Assume that there exists $\sigma \in
\Aut(R)$ lifting the $3$-rd power Frobenius automorphism of $k$.
Let $A$ be an abelian scheme of relative dimension $g$
over $R$, which is assumed to have ordinary reduction, and
let $\pol$ be an ample symmetric line bundle of degree $1$ on $A$.
We set $Z_n=(\xz/n \xz)_R^g$ for an integer $n \geq 1$.
Assume that we are given
a symmetric theta structure $\Theta_4$ of type $Z_4$ for the pair
$(A, \pol^4)$.
Let $(a_u)_{u \in Z_4}$ denote the theta null point with respect to
the theta structure $\Theta_4$.
In the following we identify $Z_2$ with its image in $Z_4$ under the
morphism which maps component-wise $1 \mapsto 2$. We define
\begin{eqnarray*}
S= \{ (x,y,z) \in Z_4^3 \,\mid\, (x-2y,x+y-z,x+y+z) \in Z_2^3 \}.
\end{eqnarray*}
For $(x_1,y_1, z_1),(x_2,y_2,z_2) \in S$ we denote $(x_1,y_1,z_1) \sim (x_2,y_2,z_2)$
if there exists a permutation matrix $P \in \mat{3}{\xz}$ such that
\[
(x_1-2y_1,x_1+y_1-z_1,x_1+y_1+z_1) = (x_2-2y_2,x_2+y_2-z_2,x_2+y_2+z_2) P.
\]
\begin{theorem}
\label{kop}
Assume that $A$ is the canonical lift of $A_k$.
For
$(x,y_1,z_1),(x,y_2,z_2) \in S$ such that 
$(x,y_1,z_1) \sim (x,y_2,z_2)$ one has
\begin{eqnarray*}
\sum_{ u \in Z_2} a_{y_1+u}^\sigma a_{z_1+u}
=\sum_{ v \in Z_2} a_{y_2+v}^\sigma a_{z_2+v}.
\end{eqnarray*}
\end{theorem}
\begin{proof}
There exists a unique theta structure $\Theta_2$
of type $Z_2$ for $(A, \pol^2)$ which is $2$-compatible
with the given theta structure $\Theta_4$ (see \cite[$\S$2,Rem.1]{mu66}).
Now assume that we have chosen an isomorphism
\begin{eqnarray}
\label{tors3}
Z_3 \stackrel{\sim}{\rightarrow} A[3]^{\mathrm{et}}.
\end{eqnarray}
In order to do so we may have to extend locally-\'{e}tale the base ring $R$.
Note that $\sigma$ admits a unique continuation to local-\'{e}tale extensions.
Our assumption is justified by the following observation. As we shall
see lateron, the resulting theta relations have coefficients in $\xz$
and hence are defined over the original ring $R$.
\newline\indent
By \cite[Th.2.2]{ca05a} the isomorphism (\ref{tors3}) determines a canonical theta
structure $\Theta_3^{\mathrm{can}}$ of type $Z_3$ for $\pol^3$.
By Lemma \ref{prodcantheta} there exist semi-canonical product theta structures
$\Theta_6=\Theta_2 \times \Theta_3^{\mathrm{can}}$ and
$\Theta_{12}=\Theta_4 \times \Theta_3^{\mathrm{can}}$
of type $Z_6$ and $Z_{12}$ for $\pol^6$ and $\pol^{12}$, respectively.
By \cite[Th.5.1]{ca05a} and Lemma \ref{dminus1} the canonical theta
structure $\Theta_3^{\mathrm{can}}$ is symmetric. We conclude by
Lemma \ref{prodcompat} that the theta structures $\Theta_2$, $\Theta_4$, $\Theta_6$ and
$\Theta_{12}$ are compatible in the sense of \cite[$\S$5.3]{ca05b}.
For the following we assume that we have chosen rigidifications
for the line bundles $\pol^i$ and theta invariant isomorphisms
\[
\mu_i: \pi_* \pol^i \stackrel{\sim}{\rightarrow} V(Z_i)= \underline{\mathrm{Hom}}(Z_i, \mathcal{O}_R),
\]
where $i=2,4,6,12$ and $\pi:A \rightarrow \spec{R}$ denotes the
structure morphism.
Our choice determines theta functions $q_{\pol^i} \in V(Z_i)$ which
interpolate the coordinates of the theta null point with respect to
$\Theta_i$ (see \cite[$\S$1]{mu66}).
Let $\{ \delta_w \}_{w \in Z_2}$ denote the Dirac basis 
of the module of finite
theta functions $V(Z_2)$.
Let now $(x_0,y_i,z_i) \in S$ where $i=1,2$ and set
\[
(a_i,b_i,c_i)=(x_0-2y_i,x_0+y_i-z_i,x_0+y_i+z_i).
\]
Suppose that $(x_0,y_1,z_1) \sim (x_0,y_2,z_2)$,
i.e. there exists
a permutation matrix $P \in \mat{3}{\xz}$ such that
\begin{eqnarray}
\label{permut}
(a_1,b_1,c_1) = (a_2,b_2,c_2) P.
\end{eqnarray}
For $i=1,2$ we set
\[
S^{(i)}_{x_0}=\{ (x,y,z) \in S \hspace{0.1cm} |  \hspace{0.1cm}
(x=x_0)
\wedge (x-2y,x+y-z,x+y+z)=(a_i,b_i,c_i) \}.
\]
By Theorem \ref{threemult}
there exists a $\lambda \in R^*$ such that
\begin{eqnarray}
\label{prodbasis}
\lefteqn{ \big( \delta_{a_i} \star \delta_{b_i} \star \delta_{c_i}
  \big) (x_0)} \\
\nonumber & & = \lambda \sum_{ (x,y,z) \in S^{(i)}_{x_0} } \delta_{a_i}(x-2y) \delta_{b_i}(x+y-z)
\delta_{c_i}(x+y+z) q_{\pol^{12}}(y) q_{\pol^4}(z) \\
\nonumber & & = \lambda \sum_{t \in Z_2}
q_{\pol^{12}}(y_i+t) q_{\pol^4}(z_i+t).
\end{eqnarray}
It follows by Theorem \ref{cantwisttheta} and Lemma \ref{descentone}
that there exists an $\alpha \in R^*$ such that
\begin{eqnarray}
\label{desc}
q_{\pol^{12}}(z)= \alpha q_{\pol^4}(z)^{\sigma}
\end{eqnarray}
for all $z \in Z_4$.
Combining the equations (\ref{prodbasis}) and (\ref{desc})
we conclude
that there exists $\lambda \in R^*$ such that
\begin{eqnarray}
\label{form}
\big( \delta_{a_i} \star \delta_{b_i} \star \delta_{c_i} \big) (x_0)
 = \lambda \sum_{t \in Z_2}
q_{\pol^4}(y_i+t)^{\sigma} q_{\pol^4}(z_i+t).
\end{eqnarray}
The commutativity of the $\star$-product and equality (\ref{permut})
imply that
\begin{eqnarray}
\label{comm}
\big( \delta_{a_1} \star \delta_{b_1} \star \delta_{c_1} \big)(x_0) =
\big( \delta_{a_2} \star \delta_{b_2} \star \delta_{c_2} \big)(x_0).
\end{eqnarray}
As a consequence of the equalities (\ref{form}) and (\ref{comm}) we have
\begin{eqnarray*}
\sum_{u \in Z_2} q_{\pol^4}(y_1+u)^{\sigma} q_{\pol^4}(z_1+u)
= \sum_{v \in Z_2} q_{\pol^4}(y_2+v)^{\sigma} q_{\pol^4}(z_2+v).
\end{eqnarray*}
This completes the proof of Theorem \ref{kop}.
\end{proof}
Note that by symmetry one has $a_u=a_{-u}$ for all $u \in Z_4$.
For the sake of completeness we also give the well-known
higher dimensional modular equations of level $4$ which generalize Riemann's
relation.
Let
\begin{eqnarray*}
S'= \{ (v,w,x,y) \in Z_4^4 \hspace{0.2cm} |  \hspace{0.2cm} (v+w,v-w,x+y,x-y) \in Z_2^4 \}.
\end{eqnarray*}
For $(v_1,w_1,x_1,y_1),(v_2,w_2,x_2,y_2) \in S'$ we write
$(v_1,w_1,x_1,y_1) \sim (v_2,w_2,x_2,y_2)$
if there exists a permutation matrix $P \in \mat{4}{\xz}$ such that
\[
(v_1+w_1,v_1-w_1,x_1+y_1,x_1-y_1)=(v_2+w_2,v_2-w_2,x_2+y_2,x_2-y_2)P.
\] 
\begin{theorem}
\label{riemann}
For $(v_1,w_1,x_1,y_1),(v_2,w_2,x_2,y_2) \in S'$ such that 
$(v_1,w_1,x_1,y_1) \sim (v_2,w_2,x_2,y_2)$, 
the following equality holds
\begin{eqnarray*}
\sum_{t \in Z_2} a_{v_1+t} a_{w_1+t} \sum_{s \in Z_2} a_{x_1+s} a_{y_1+s}
=\sum_{t \in Z_2} a_{v_2+t} a_{w_2+t} \sum_{s \in Z_2} a_{x_2+s} a_{y_2+s}.
\end{eqnarray*}
\end{theorem}
A proof of the above theorem can be found in \cite[$\S$3]{mu66}.

\subsubsection{Theta null values in dimensions $1$ and $2$}
\label{equations}

In this section we make the equations of Theorem~\ref{kop} and 
Theorem~\ref{riemann} explicit in the case of dimensions $1$ and $2$.
Let $\xf_q$ be a finite field of characteristic $3$ having $q$ elements 
and let $R=W(\xf_q)$ denote the Witt vectors with values in $\xf_q$.
There exists a canonical lift $\sigma \in \mathrm{Aut}(R)$ of the $3$-rd 
power Frobenius of $\xf_q$.
Let $A$ be an abelian scheme over $R$ with ample symmetric line bundle 
$\pol$ of degree $1$ on $A$.

\medskip

\noindent{\bf Dimension 1.}
Suppose that $A$ is a proper smooth elliptic curve over $R$, and let 
$(a_0:a_1:a_2:a_3)$ be the theta null point with respect to a 
symmetric theta structure of type $(\xz/4\xz)_R$ for $(A, \pol^4)$
where $\pol=\pol(0_A)$ and $0_A$ denotes the zero section of $A$.
By symmetry we have $a_1 = a_3$, and Theorem~\ref{riemann} implies that  
the projective point $(a_0:a_1:a_2)$ lies on the smooth genus~$3$ curve 
$\mathcal{A}_1(\Theta_4) \subseteq
\proj{\xz[\frac{1}{2},x_0,x_1,x_2]}=\xp^2_{\xz[\frac{1}{2}]}$
with defining equation 
\begin{equation}
\label{eqn_riemann_g1}
(x_0^2+x_2^2)x_0x_2 = 2x_1^4.
\end{equation}
The latter classical equation is known as \emph{Riemann's relation}.
We remark that the points on $\mathcal{A}_1(\Theta_4)$ give the moduli of
elliptic curves with symmetric $4$-theta structure.
\newline\indent
Now assume that $A$ has ordinary reduction and that $A$ is the
canonical lift of $A_{\xf_q}$.
Theorem~\ref{kop} implies that the coordinates of the
projective point $(a_0:a_1:a_2)$ satisfy
the equation
\begin{equation}
\label{eqn_corresp_g1}
x_0 y_2 + x_2 y_0 = 2x_1 y_1,
\end{equation}
where $x_i=a_i$ and $y_i=a_i^{\sigma}$ for $i=0,1,2$.

\medskip

\noindent{\bf Dimension 2.}
Now suppose that $A$ has relative dimension $2$ over $R$ and that we
are given a symmetric 
theta structure of type $(\xz/4\xz)^2$ for the pair $(A, \pol^4)$.
Let $(a_{ij})_{(i,j) \in (\xz/4\xz)^2}$ denote the theta null point 
with respect to the latter theta structure. By symmetry we have
\begin{eqnarray*}
a_{11}=a_{33}, \hspace{0.2cm}
a_{10}=a_{30}, \hspace{0.2cm}
a_{01}=a_{03}, \hspace{0.2cm}
a_{13}=a_{31}, \hspace{0.2cm}
a_{32}=a_{12}, \hspace{0.2cm}
a_{21}=a_{23}.
\end{eqnarray*}
The $2$-dimensional analogue of \emph{Riemann's equation}
(\ref{eqn_corresp_g1}) are the equations
\begin{equation}
\label{RiemEqns}
\begin{array}{c}
(x_{00}^2+x_{02}^2+x_{20}^2+x_{22}^2)(x_{00}x_{02}+x_{20}x_{22}) 
  = 2(x_{01}^2+x_{21}^2)^2 \\
(x_{00}^2+x_{02}^2+x_{20}^2+x_{22}^2)(x_{00}x_{20}+x_{02}x_{22}) 
  = 2(x_{10}^2+x_{12}^2)^2 \\
(x_{00}^2+x_{02}^2+x_{20}^2+x_{22}^2)(x_{00}x_{22}+x_{20}x_{02}) 
  = 2(x_{11}^2+x_{13}^2)^2 \\
(x_{00}x_{20}+x_{02}x_{22})(x_{00}x_{22}+x_{02}x_{20}) 
  = 4x_{01}^2 x_{21}^2 \\
(x_{00}x_{02}+x_{20}x_{22})(x_{00}x_{22}+x_{02}x_{20}) 
  = 4x_{10}^2 x_{12}^2 \\
(x_{00}x_{02}+x_{20}x_{22})(x_{00}x_{20}+x_{02}x_{22}) 
  = 4x_{11}^2 x_{13}^2 \\
(x_{00}^2+x_{02}^2+x_{20}^2+x_{22}^2)x_{13}x_{11}
  = (x_{12}^2+x_{10}^2)(x_{01}^2+x_{21}^2) \\
(x_{00}^2+x_{02}^2+x_{20}^2+x_{22}^2)x_{01}x_{21} 
  = (x_{12}^2+x_{10}^2)(x_{11}^2+x_{13}^2) \\
(x_{00}^2+x_{02}^2+x_{20}^2+x_{22}^2)x_{10}x_{12}
  = (x_{01}^2+x_{21}^2)(x_{11}^2+x_{13}^2) \\
(x_{02}x_{20}+x_{00}x_{22})x_{11}x_{13}
  = 2 x_{01}x_{10}x_{21}x_{12} \\
(x_{20}x_{00}+x_{22}x_{02})x_{10}x_{12} 
  = 2x_{11}x_{13}x_{21}x_{01} \\
(x_{00}x_{02}+x_{20}x_{22})x_{21}x_{01} 
  = 2x_{11}x_{13}x_{10}x_{12} \\
(x_{02}x_{20}+x_{00}x_{22})(x_{01}^2+x_{21}^2)
  = 2 x_{10}x_{12}(x_{11}^2+x_{13}^2) \\
(x_{00}x_{02}+x_{20}x_{22})(x_{11}^2+x_{13}^2)
  = 2x_{10}x_{12}(x_{01}^2+x_{21}^2) \\
(x_{02}x_{20}+x_{00}x_{22})(x_{10}^2+x_{12}^2)
  = 2 x_{21}x_{01}(x_{11}^2+x_{13}^2) \\
(x_{20}x_{00}+x_{22}x_{02})(x_{13}^2+x_{11}^2)
  = 2x_{21}x_{01}(x_{10}^2+x_{12}^2) \\
(x_{20}x_{00}+x_{22}x_{02})(x_{21}^2+x_{01}^2)
  = 2x_{11}x_{13}(x_{10}^2+x_{12}^2) \\
(x_{00}x_{02}+x_{20}x_{22})(x_{12}^2+x_{10}^2)
  = 2x_{11}x_{13}(x_{01}^2+x_{21}^2) \\
x_{01}x_{21}(x_{01}^2+x_{21}^2) = x_{10}x_{12}(x_{10}^2+x_{12}^2) \\
x_{01}x_{21}(x_{01}^2+x_{21}^2) = x_{11}x_{13}(x_{11}^2+x_{13}^2).
\end{array}
\end{equation}
By Theorem \ref{riemann} the point $(a_{ij})_{(i,j) \in (\xz/4\xz)^2}$
is a solution of the equations (\ref{RiemEqns}),
i.e. the above equations hold for $x_{ij}=a_{ij}$.
The latter equations determine a three dimensional subscheme $\mathcal{A}_2(\Theta_4)$
of the projective space
\[
\xp^{9}_{\xz[\frac{1}{2}]}=\mathrm{Proj} \big(
\xz[\frac{1}{2},x_{00},x_{01},x_{02},x_{10},x_{11},x_{12},x_{13},x_{20},x_{21},x_{22}]
\big).
\]
The points on $\mathcal{A}_2(\Theta_4)$ give the moduli
of abelian surfaces with symmetric theta structure of type 
$(\xz/4\xz)^2$.
We remark that the point
\[
(a_{00}:a_{01}:a_{02}:a_{10}:a_{11}:a_{12}:a_{13}:a_{20}:a_{21}:a_{22})
\in \xp^{9}_{\xz[\frac{1}{2}]}(R)
\]
is a solution of the 
equations (\ref{RiemEqns}) if and only if the projective coordinates
$$
\begin{array}{cc}
& \big(  a_{00}^2 + a_{02}^2 + a_{20}^2 + a_{22}^2 : 
   2(a_{01}^2 + a_{21}^2) : 
   2(a_{12}^2 + a_{10}^2) : 
   2(a_{11}^2 + a_{13}^2) \big), \\
& \big(
   a_{01}^2 + a_{21}^2 : 
   a_{00}a_{02} + a_{20}a_{22} : 
   2a_{11}a_{13} : 
   2a_{10}a_{12} \big), \\
& \big(
   a_{12}^2 + a_{10}^2 : 
   2a_{11}a_{13} : 
   a_{00}a_{20} + a_{02}a_{22} : 
   2a_{01}a_{21} \big), \\
& \big(
    a_{11}^2 + a_{13}^2 : 
    2a_{10}a_{12} : 
    2a_{01}a_{21} : 
    a_{00}a_{22} + a_{02}a_{20}\big)
\end{array}
$$
describe the same point in $\xp^3_{\xz[\frac{1}{2}]}(R)$.
In fact the above formulas define a morphism to the space of
abelian surfaces with $2$-theta structure which embeds in $\xp^3_{\xz[\frac{1}{2}]}$.
Together with the Riemann equations, the following corollary of 
Theorem~\ref{kop} forms the basis of our construction algorithm for 
CM abelian surfaces.
\begin{corollary}
\label{genus2}
Assume that $A$ has ordinary reduction and that $A$ is the canonical lift 
of $A_{\xf_q}$. Let $(a_{ij})$ denote the theta null point of $A$
with respect to a given symmetric $4$-theta structure. Then the coordinates of the point
\[
(a_{00}:a_{01}:a_{02}:a_{10}:a_{11}:a_{12}:a_{13}:a_{20}:a_{21}:a_{22})
\in \xp^{9}_{\xz[\frac{1}{2}]}(R)
\]
satisfy
the following relations
\begin{equation}
\label{CorrEqns}
\begin{array}{c}
x_{00} y_{02} + x_{02} y_{00} + x_{20} y_{22} + x_{22} y_{20} 
  - 2(x_{01} y_{01} + x_{21} y_{21}) =0 \\ 
x_{00} y_{20} + x_{20} y_{00} + x_{02} y_{22} + x_{22} y_{02} 
  - 2(x_{10} y_{10} + x_{12} y_{12}) =0 \\ 
x_{00} y_{22} + x_{22} y_{00} + x_{02} y_{20} + x_{20} y_{02} 
  - 2(x_{13} y_{13} + x_{11} y_{11}) =0 \\ 
x_{01} y_{21} + x_{21} y_{01} - (x_{12} y_{10} + x_{10} y_{12})=0 \\ 
x_{01} y_{21} + x_{21} y_{01} - ( x_{11} y_{13} + x_{13} y_{11})=0,
\end{array}
\end{equation}
where $x_{ij} = a_{ij}$ and $y_{ij} = a_{ij}^\sigma$.
\end{corollary}

\subsection{Galois properties of the canonical theta structure}
\label{galcantheta}

For our notation and standard definitions we refer to \cite{mu66}
and \cite{ca05a}.
Let $R$ be a complete noetherian local ring with
perfect residue field $k$ of characteristic $p>2$.
Suppose that we are given an abelian scheme $A$ over $R$ which has ordinary reduction.
Let $\pol$ be an ample symmetric line bundle of degree $1$ on $A$.
We set $q=p^d$ where $d \geq 1$ is an integer.
Assume that there exists a $\sigma \in \mathrm{Aut}(R)$ lifting the $q$-th power
Frobenius automorphism of $k$.
Recall that there exists a canonical lift $F:A \rightarrow A^{(q)}$ of
the relative $q$-Frobenius morphism and a canonical ample symmetric line bundle
$\pol^{(q)}$ of degree $1$ on $A^{(q)}$ such that $F^* \pol^{(q)} \cong \pol^{q}$
(see \cite[$\S$5]{ca05a}).
Let $A^{(\sigma)}$ be defined by the Cartesian diagram
\begin{eqnarray}
\label{defpullback}
\xymatrix{
A^{(\sigma)} \ar@{->}[r]^{\mathrm{pr}} \ar@{->}[d] & A \ar@{->}[d] \\
\spec{R} \ar@{->}[r]^{\spec{\sigma}} & \spec{R},}
\end{eqnarray}
where the right hand vertical arrow is the structure morphism.
Let $\pol^{(\sigma)}$ be the pull back of $\pol$ along the morphism
$\mathrm{pr}:A^{(\sigma)} \rightarrow A$ which is defined by
the diagram (\ref{defpullback}).
\newline\indent
Now let $n \geq 1$ be a natural number such that $(n,p)=1$, i.e. the
numbers $n$ and $p$ are coprime.
Assume that we are given a symmetric theta structure $\Theta_n$ of
type $Z_n = (\xz / n \xz)_R^g$
for $\pol^n$ where $g=\mathrm{dim}_R(A)$.
We denote by $\pol_n^{(\sigma)}$ the $n$-th power
of $\pol^{(\sigma)}$.
We obtain a theta structure $\Theta_n^{(\sigma)}$ for $\pol_n^{(\sigma)}$
on $A^{(\sigma)}$ by extension of scalars along $\spec{\sigma}$
applied to the theta structure $\Theta_n:G(Z_n) \stackrel{\sim}{\rightarrow}
G(\pol_n)$ and by chaining with the natural
isomorphism $Z_n \stackrel{\sim}{\rightarrow} Z_{n,\sigma}$ and
the inverse of its dual.
Assume that $A$ is the canonical lift of $A_k$.
Our assumption implies that the abelian schemes $A^{(q)}$ and $A^{(\sigma)}$ are canonical lifts
which are canonically isomorphic over the residue field $k$.
In the following we will assume that the special fibers of $A^{(q)}$
and $A^{(\sigma)}$ are indeed equal.
There exists a canonical isomorphism
$\tau: A^{(q)} \stackrel{\sim}{\rightarrow} A^{(\sigma)}$ over $R$
lifting the identity on special fibers.
\newline\indent
We claim that $\tau^* \pol^{(\sigma)} \cong \pol^{(q)}$.
We set $\mathcal{M}=\tau^* \pol^{(\sigma)} \otimes (
\pol^{(q)})^{-1}$.
It follows by the definition of $\pol^{(q)}$
that the class of $\mathcal{M}$ reduces
to the trivial class.
Note that $\tau^* \pol^{(\sigma)}$ is symmetric.
By \cite[Th.5.1]{ca05a} also the line bundle $\pol^{(q)}$ is symmetric.
As a consequence the line bundle $\mathcal{M}$ is symmetric
and gives an element of
$\mathrm{Pic}^0_{A^{(q)}/R}[2](R)$.
We observe that the group $\mathrm{Pic}^0_{A^{(q)}/R}[2]$ is finite \'{e}tale
because of the assumption $p>2$.
We conclude by the connectedness of the ring $R$
that the class of $\mathcal{M}$ is the trivial class. Hence our claim follows.
\newline\indent
By the above discussion there exists an isomorphism
$\gamma: \tau^* \pol^{(\sigma)} \stackrel{\sim}{\rightarrow} \pol^{(q)}$.
We define a $\xg_{m,R}$-invariant morphism of theta groups $\tau^*:G(\pol^{(\sigma)}) \rightarrow
G(\pol^{(q)})$ by setting
$(x, \varphi) \mapsto \big( y, T_y^* \gamma \circ \tau^* \varphi \circ \gamma^{-1} \big)$
where $y = \tau^{-1}(x)$ and $\varphi: \pol^{(\sigma)}
\stackrel{\sim}{\rightarrow} T_x  \pol^{(\sigma)}$.
Obviously, our definition is independent of the choice of
$\gamma$. For trivial reasons the morphism $\tau^*$ gives an isomorphism.
\begin{theorem}
\label{cantwisttheta}
There exists a canonical theta structure $\Theta_n^{(q)}$ of type
$Z_n$ for $\pol_n^{(q)}$ depending on $\Theta_n$ such that
\begin{eqnarray}
\label{pullbackprop}
\tau^* \circ \Theta_n^{(\sigma)} = \Theta_n^{(q)}.
\end{eqnarray}
\end{theorem}
\begin{proof}
Assume that we have chosen an isomorphism
\begin{eqnarray}
\label{trivq}
Z_q = (\xz / q \xz)^g_R \stackrel{\sim}{\rightarrow} A[q]^{\mathrm{et}},
\end{eqnarray}
where $A[q]^{\mathrm{et}}$ denotes the maximal \'{e}tale
quotient of $A[q]$.
In order to do so we may have to extend
$R$ locally-\'{e}tale.
By \cite[Th.2.2]{ca05a} there exists a canonical theta structure $\Theta^{\mathrm{can}}_q$
of type $Z_q$ for the pair $\big( A, \pol^q \big)$
depending on the isomorphism \rm{(\ref{trivq})}.
We remark that the canonical theta structure is symmetric
by \cite[Th.5.1]{ca05a} and Lemma \ref{dminus1}.
By Lemma \ref{prodcantheta} there exists a semi-canonical symmetric product theta structure
$\Theta_{nq}=\Theta_n \times \Theta^{\mathrm{can}}_q$
of type $Z_{nq}=(\xz / nq \xz)^g_R$
for the pair $\big( A,\pol_{nq} \big)$ where $\pol_{nq}=\pol^{\otimes nq}$.
It follows from \cite[Prop.5.3]{ca05b} that the theta structure
$\Theta_{nq}$ descends along the Frobenius lift $F$ to a canonical theta
structure $\Theta_n^{(q)}$ for $\pol_n^{(q)}$.
We choose $\Theta_n^{(q)}=\Theta_{nq}(\mathrm{id})$
where our notation is as in \cite[$\S$5.2]{ca05b}. 
In the following we will prove that $\Theta_n^{(q)}$ has the desired
pull back property (\ref{pullbackprop}).
\newline\indent
First we check the pull back property for the induced
Lagrangian structures.
Let $\delta_n$, $\delta_n^{(q)}$ and $\delta_n^{(\sigma)}$ be the Lagrangian
structures which are induced by $\Theta_n$, $\Theta_n^{(q)}$ and
$\Theta_n^{(\sigma)}$, respectively.
We claim that
\begin{eqnarray}
\label{pullbacklag}
\delta_n^{(\sigma)} = \tau \circ \delta_n^{(q)}.
\end{eqnarray}
As $\Theta_n^{(q)}=\Theta_{nq}(\mathrm{id})$
(notation as in \cite[$\S$5.2]{ca05b})
the restriction of the Lagrangian structure $\delta_n^{(q)}$ to $Z_n$
equals the restriction of $F \circ \delta_n$. As a consequence the
restrictions of the morphisms
$\tau \circ \delta_n^{(q)}$ and $\delta_n^{(\sigma)}$
coincide on the special fiber.
By general theory the reduction functor on the category of
finite \'{e}tale schemes over
$R$ gives an equivalence of categories.
We conclude that the equality (\ref{pullbacklag})
restricted to $Z_n$ is true over $R$.
The equality for $Z_n^D$ can be proven analogously. Note that $Z_n^D$ is
\'{e}tale because of the assumption $(n,p)=1$.
Hence the claim follows.
\newline\indent
It remains to show, on top of Lagrangian structures, the equality of theta structures as claimed
in (\ref{pullbackprop}).
By \cite[Prop.4.5]{ca05a}
the theta structures $\Theta_n$, $\Theta_n^{(\sigma)}$ and $\Theta_n^{(q)}$ give
rise to sections $s$, $s^{(\sigma)}$ and $s^{(q)}$ of theta exact
sequences
\begin{eqnarray}
\label{thetacommute}
\xymatrix{
0 \ar@{->}[r] & \xg_{m,R} \ar@{->}[dd]^{\mathrm{id}} \ar@{->}[r] & G(\pol^{(q)}_n)
\ar@{->}[rr] \ar@{<-}[dd]_{\tau^*} & & H(\pol^{(q)}_n) \ar@{->}[dd]^{\tau}
\ar@{->}[r] & 0 \\
& & & Z_n \ar@{->}[ul]^{s^{(q)}} \ar@{->}[ur]^{\delta_n^{(q)}}
\ar@{->}[dl]^{s^{(\sigma)}} \ar@{->}[dr]^{\delta_n^{(\sigma)}} & & \\
0 \ar@{->}[r] & \xg_{m,R} \ar@{->}[r] & G(\pol^{(\sigma)}_n)
\ar@{->}[rr] & & H(\pol^{(\sigma)}_n) \ar@{->}[r] & 0 .}
\end{eqnarray}
We claim that the diagram (\ref{thetacommute}) commutes.
By definition, the squares of the above diagram commute.
By equation (\ref{pullbacklag}) the right hand triangle commutes.
It remains to show that
\begin{eqnarray}
\label{pullbacksect}
\tau^* \circ s^{(\sigma)} = s^{(q)}.
\end{eqnarray}
The difference of
$\tau^* \circ s^{(\sigma)}$ and $s^{(q)}$ gives a point
$\varphi \in \mathrm{Hom}_R(Z_n,\xg_{mR})=\mu_{n,R}(R)$.
It suffices to show that the point $\varphi$ reduces to the neutral
element of $\mu_{n,k}$, because
the group $\mu_{n,R}$ is \'{e}tale and the ring $R$ is connected.
In the following we prove that $\varphi$ has trivial reduction.
Consider the diagram
\begin{eqnarray*}
\xymatrix{
G(\pol_n) \ar@{->}[r]^{\epsilon_q} \ar@{<-}[d]^{\Theta_n} &
G(\pol_{nq}) \ar@{<-}[d]^{\Theta_{nq}}  &
G(\pol_{nq})^* \ar@{->}[r] \ar@{_{(}->}[l] & G(\pol_{nq})^*/\tilde{K} \ar@{<-}[d]_{\cong} \ar@{->}[r]^{\mathrm{can},F} &
G(\pol_n^{(q)}) \ar@{<-}[d]_{\Theta_n^{(q)}} \\
G(Z_n) \ar@{->}[r]^{E_q} & G(Z_{nq}) & 
\xg_{m,R} \times Z_n \times
Z_{nq}^D \ar@{->}[r]  \ar@{->}[u]^{\cong} \ar@{_{(}->}[l] & \xg_{m,R} \times Z_n \times
(Z_{nq}^D /K) \ar@{->}[r] & G(Z_n)
}
\end{eqnarray*}
where $E_q$ and $\epsilon_q$ are defined as in \cite[$\S$5.3]{ca05b},
$K=\mathrm{Ker}(F)$ and $\tilde{K}$
is a canonical lift of $K$ to the theta group $G(\pol_{nq})$. The group
$G(\pol_{nq})^*$ is defined as the centralizer of $\tilde{K}$ in $G(\pol_{nq})$.
By Lemma \ref{prodcompat} the left hand square of the above diagram is commutative.
The lift $\tilde{K}$ is induced
by some isomorphism $\alpha:F^* \pol_n^{(q)} \stackrel{\sim}{\rightarrow}
\pol_{nq}$. Let $x \in Z_n$, $y= \delta_n(x,1)$ and $z=\delta_n^{(q)}(x,1)$.
Suppose that $s(x)=\Theta_n(1,x,1)=(y,\psi_y)$ and
$s^{(q)}(x)=\Theta_n^{(q)}(1,x,1)=(z, \gamma_z)$.
Note that $z=F(y)$.
It follows by the commutativity of the above diagram that
\begin{eqnarray}
\label{keyeq}
\psi_y^{\otimes q} =
T_y^* \alpha \circ F^* \gamma_z \circ \alpha^{-1}.
\end{eqnarray}
Equation (\ref{keyeq}) says in down-to-earth terms that,
on the special fiber, the isomorphism $\gamma_z$ is
the pull back of $\psi_y$ under the isomorphism
$\mathrm{pr}:A^{(\sigma)} \rightarrow A$ where
the latter is defined by the diagram (\ref{defpullback}).
This proves that the above character $\varphi$ is trivial on the
special fibre.
This proves the equality (\ref{pullbacksect}) and hence the diagram
(\ref{thetacommute}) is commutative.
The proof for $Z_n^D$ is analogous. 
\newline\indent
We remark that the equality (\ref{pullbackprop}) implies by means of descent
that $\Theta_n^{(q)}$ is defined over $R$.
This completes the proof of the theorem.
\end{proof}

\section{On the theory of algebraic theta functions}
\label{basic}

In this section we prove some basic facts about algebraic theta
functions which are needed in the proof of Theorem \ref{kop}.
These results are absent
from the literature. For an
introduction to algebraic theta functions we refer to \cite{mu66}.

\subsection{Symmetric theta structures}
\label{thetasym}

In this section we recall the notion of a \emph{symmetric theta structure}.
The symmetry turns out to be an essential ingredient in the proof
of the theta relations of Theorem \ref{kop}.
We give a characterization of the symmetry of a theta structure
in terms of the symmetry of the associated line bundles. This characterization is not
obvious from the definitions given in \cite[$\S$2]{mu66}.
The results of this section imply that the canonical theta
structure, whose existence is proven in \cite{ca05a},
is a symmetric theta structure. Note that our definition of
symmetry is weaker than the one given in \cite[$\S$2]{mu66}.
\newline\indent
Let $A$ be an abelian scheme over a ring $R$
and let $\pol$ be a line bundle on $A$.
Consider the morphism
\[
\varphi_{\pol}:A \rightarrow \mathrm{Pic}^0_{A/R},
\hsp x \mapsto \langle T_x^* \pol \otimes \pol^{-1} \rangle
\]
where $\langle \cdot \rangle$ denotes the class in
$\mathrm{Pic}^0_{A/R}$.
We denote the kernel of $\varphi_{\pol}$ by $A[\pol]$.
The line bundle $\pol$ is called \emph{symmetric}
if $[-1]^* \pol \cong \pol$.
\newline\indent
Now assume that we are given an isomorphism
$\psi:\pol \stackrel{\sim}{\rightarrow} [-1]^* \pol$.
We denote the theta group of the line bundle $\pol$ by $G(\pol)$.
Let $(x,\varphi) \in G(\pol)$, where $x \in A[\pol]$ and $\varphi:\pol
\stackrel{\sim}{\rightarrow} T_x^* \pol$ is an isomorphism, and let $\tau_{\varphi}$ denote the composed
isomorphism
\[
\pol \stackrel{\psi}{\rightarrow} [-1]^* \pol
\stackrel{[-1]^* \varphi}{\longrightarrow} [-1]^* T_x^* \pol =
T^*_{-x} [-1]^* \pol \stackrel{T^*_{-x} \psi^{-1}}{\longrightarrow}
T^*_{-x} \pol.
\]
One defines a morphism $\delta_{-1}:G(\pol) \rightarrow G(\pol)$
by setting $\delta_{-1}(x,\varphi)=(-x,\tau_{\varphi})$.
We remark that the definition of $\tau_{\varphi}$
does not depend on the choice of the isomorphism $\psi$.
Obviously $\delta_{-1}$ is an automorphism of order $2$ of the group $G(\pol)$.
\newline\indent
Let $K$ be a finite constant group over $R$.
We define an automorphism $D_{-1}$ of the standard theta group
$G(K)=\xg_{m,R} \times K \times K^D$ by mapping
$(\alpha,x,l) \mapsto (\alpha,-x,l^{-1})$.
Assume now that we are given a theta structure  $\Theta:G(K)
\stackrel{\sim}{\rightarrow} G(\pol)$.
\begin{definition}
The theta structure $\Theta$ is called \emph{symmetric} if the following
equality holds
\begin{eqnarray}
\label{sym}
\Theta \circ D_{-1} = \delta_{-1} \circ \Theta.
\end{eqnarray}
\end{definition}
Note that we do not assume that the line bundle $\pol$
is totally symmetric as it is done in \cite[$\S$2]{mu66}. 
In the following we will give a necessary and sufficient condition for
a theta structure to be symmetric.
Recall (see \cite[$\S$4]{ca05a}) that 
the theta structure $\Theta$ corresponds to a Lagrangian
structure of type $K$ and isomorphisms
\begin{eqnarray*}
\alpha_K:I_K^* \mathcal{M}_K \stackrel{\sim}{\rightarrow} \pol
\quad \mbox{and} \quad 
\alpha_{K^D}:I_{K^D}^* \mathcal{M}_{K^D} \stackrel{\sim}{\rightarrow} \pol,
\end{eqnarray*}
where $I_K:A \rightarrow A_K$ and $I_{K^D}:A \rightarrow A_{K^D}$
are isogenies with kernel $K$ and $K^D$,
and $\mathcal{M}_K$ and $\mathcal{M}_{K^D}$
are line bundles on $A_K$ and $A_{K^D}$, respectively.
\begin{lemma}
\label{dminus1}
The theta structure $\Theta$ is symmetric if and only if
the line bundles $\mathcal{M}_K$ and $\mathcal{M}_{K^D}$
are symmetric.
\end{lemma}
\begin{proof}
We prove that the equality (\ref{sym}) holds on the image of the
morphism
\[
s_K:K \rightarrow G(K), x \mapsto (1,x,1)
\]
if and only if the line bundle $\mathcal{M}_K$
is symmetric.
An analogous proof exists for the dual construction.
Consider the following diagram
\begin{eqnarray*}
\xymatrix{
A[\pol] \ar@{<-}[r]^{\mathrm{proj}} \ar@{<-}[d]^{j} & G(\pol) \ar@{->}[r]^{\delta_{-1}} \ar@{<-}[d]^{\Theta} &
G(\pol) \ar@{<-}[d]^{\Theta} \\
K \ar@{->}[r]^{s_K} & G(K) \ar@{->}[r]^{D_{-1}} & G(K).
}
\end{eqnarray*}
Here the morphism $j$ denotes the inclusion induced
by the Lagrangian structure, which is part of the theta structure $\Theta$.
By \cite[Prop.4.2 and Prop.4.5]{ca05a} we have
$\Theta(1,x,1)= \big( j(x), T_{j(x)}^* \alpha_K \circ \alpha_K^{-1} \big)$.
We conclude that
\[
\delta_{-1} \big( \Theta (1,x,1) \big) = \left( -j(x),  T_{-j(x)}^* \psi^{-1} \circ
[-1]^* ( T_{j(x)}^* \alpha_K \circ \alpha_K^{-1} ) \circ \psi \right)
\]
where $\psi:\pol \stackrel{\sim}{\rightarrow} [-1]^* \pol$ is an
isomorphism as above.
On the other hand one has
\[
\Theta \big( D_{-1} (1,x,1) \big) = \big( -j(x), T_{-j(x)}^* \alpha_K \circ \alpha_K^{-1} \big).
\]
Hence the equation (\ref{sym}) restricted to elements of the form
$(1,x,1)$ translates as
\[
T_{-j(x)}^* \psi^{-1} \circ
[-1]^* ( T_{j(x)}^* \alpha_K \circ \alpha_K^{-1} ) \circ \psi= T_{-j(x)}^* \alpha_K \circ \alpha_K^{-1}.
\]
The latter equality is equivalent to
\[
[-1]^* \alpha_K^{-1} \circ \psi \circ \alpha_K = T_{-j(x)}^* ( [-1]^* \alpha_K^{-1} \circ \psi \circ \alpha_K).
\]
The latter equality means that the composed isomorphism
\[
I_K^* \mathcal{M}_K \stackrel{\alpha_K}{\rightarrow} \pol  \stackrel{\psi}{\rightarrow} [-1]^* \pol
\stackrel{[-1]^*\alpha_K^{-1}}{\longrightarrow} [-1]^* I_K^* \mathcal{M}_K = I_K^* [-1]^* \mathcal{M}_K
\]
is invariant under $T_{-j(x)}^*$ for all $x \in K$.
This is true if and only if this isomorphism equals the pull back
of an isomorphism $\mathcal{M}_K \stackrel{\sim}{\rightarrow} [-1]^*
\mathcal{M}_K$ along $I_K$.
Thus the lemma is proven.
\end{proof}

\subsection{Product theta structures}
\label{thetanm}

The construction of \emph{product theta structures} is considered as known to the experts.
But the reader should be aware of the fact that
a product theta structure of given theta structures does not always exist.
In this section we clarify the situation by proving
the existence of a product theta structure
under a reasonable coprimality assumption.
We provide detailed proofs because of the lack of a suitable reference.
\newline\indent
Let $A$ be an abelian scheme of relative dimension $g$ over a ring
$R$ and let $\pol$ be an ample symmetric line bundle of degree $1$ on $A$.
For an integer $n \geq 1$ we set $Z_{n}=(\xz / n \xz)_R^g$. 
Now let $n,m \geq 1$ be integers
such that $(n,m)=1$, i.e. the numbers $n$ and $m$ are coprime.
Assume we are given theta structures
\[
\Theta_n:G(Z_{n}) \stackrel{\sim}{\rightarrow} G(\pol^n) \quad
\mbox{and} \quad
\Theta_m:G(Z_m) \stackrel{\sim}{\rightarrow} G(\pol^m).
\]
We consider $Z_n$ and $Z_m$ as subgroups of $Z_{nm}$ via the
morphisms that map component-wise $1 \mapsto m$ and $1 \mapsto n$, respectively. 
\begin{lemma}
\label{prodcantheta}
There exists a natural product theta structure
\[
\Theta_{nm}:G(Z_{nm}) \stackrel{\sim}{\rightarrow} G(\pol^{nm})
\]
depending on the theta structures $\Theta_n$ and $\Theta_m$.
\end{lemma}
\begin{proof}
Let $\epsilon_n$, $\epsilon_m$,
$E_n$, $E_m$, $\eta_n$, $\eta_m$, $H_n$ and $H_m$ be defined as in \cite[$\S$5.3]{ca05b}.
We claim that for all $g \in G(\pol^n)$ and $h \in G(\pol^m)$ we have
\begin{eqnarray}
\label{commutat}
\epsilon_n(h) \epsilon_m(g)=\epsilon_m(g) \epsilon_n(h),
\end{eqnarray}
where the product is taken in $G(\pol^{nm})$.
Let $\delta_n$ and $\delta_m$ denote the Lagrangian structures that
are induced by $\Theta_n$ and $\Theta_m$. We set $\delta_{nm}=\delta_n
\times \delta_m$.
Condition (\ref{commutat}) is equivalent to
\[
\mathrm{e}_{\pol^{nm}} \big( \delta_{nm}(x_g,l_g),\delta_{nm}(x_h,l_h) \big)=1,
\]
where the elements $\delta_{nm}(x_g,l_g)$ and $\delta_{nm}(x_h,l_h)$ are the images of $\epsilon_m(g)$
and $\epsilon_n(h)$, respectively, under the natural projection $G(\pol^{nm})
\rightarrow H(\pol^{nm})$.
The vanishing of the commutator pairing follows from the bilinearity
and the assumption $(n,m)=1$.
This proves the above claim.
As a consequence there exists a canonical morphism of groups
\[
\epsilon:G(\pol^n) \times
G(\pol^m) \rightarrow G(\pol^{nm}) \quad \mbox{given by} \quad
(g,h) \mapsto \epsilon_m(g) \epsilon_n(h).
\]
Because of our assumption $(n,m)=1$ the subgroup
$C=\mathrm{ker}(\epsilon)$ is contained in the subtorus $\xg_{m,S}
\times \xg_{m,S}$ of $G(\pol^n) \times G(\pol^m)$.
In the following we will prove that $\epsilon$ is surjective.
Consider the diagram
\begin{eqnarray*}
\xymatrix{ G(\pol^n) \times G(\pol^m) \ar@{->}[rr]^{\epsilon}
  \ar@{->}[d]_{\pi_n \times \pi_m} & &
  G(\pol^{nm}) \ar@{->}[d]^{\pi} \\
  H(\pol^n) \times H(\pol^m)
  \ar@{->}[rr]^{\mathrm{can}} \ar@/_1pc/[u]_s & &
  H(\pol^{nm}) }
\end{eqnarray*}
where $\pi$ and $\pi_n \times \pi_m$ denote the natural
projections. Let
$s$ be the section of $\pi_n \times \pi_m$ induced by the theta structures
$\Theta_n$ and $\Theta_m$.
We have $\pi \circ \epsilon = \pi_ n \times \pi_m$ (up to
canonical isomorphism).
As a consequence we have $\pi \circ \epsilon \circ s
= (\pi_n \times \pi_m) \circ s = \mathrm{id}$.
We conclude that
$\pi \circ \epsilon \circ s \circ \pi = \pi$.
Let $g \in G(\pol^{nm})$. Then by the latter equality the group
element $\epsilon(s(\pi(g)))$ differs from $g$ by a
unit. Hence the morphism $\epsilon$ maps a suitable multiple
of $s( \pi(g) )$ to $g$. This implies the
surjectivity of $\epsilon$.
As a consequence $\epsilon$ induces an isomorphism
\[
\tilde{\epsilon}:\big( G(\pol^n) \times G(\pol^m) \big) / C \rightarrow G(\pol^{nm})
\]
By the same reasoning as above one can define a natural morphism
$E:G(Z_n) \times G(Z_m) \rightarrow G(Z_{nm})$,
and it is readily verified that the induced morphism
$\tilde{E}:\big( G(Z_n) \times G(Z_m) \big) / C \rightarrow G(Z_{nm})$
is an isomorphism of groups.
Let $\Theta_{nm}$ denote the composed isomorphism 
\[
G(Z_{nm}) \stackrel{\tilde{E}^{-1}}{\longrightarrow} \big( G(Z_n)
\times G(Z_m) \big) / C \stackrel{\Theta_n \times \Theta_m}{\longrightarrow}
\big( G(\pol^n) \times G(\pol^m) \big) / C \stackrel{\tilde{\epsilon}}{\rightarrow} G(\pol^{nm}).
\]
The morphism $\Theta_{nm}$ establishes the theta structure whose
existence is claimed in the lemma.
\end{proof}
Now let $\Theta_{nm}$ be as in Lemma
\ref{prodcantheta} and define the $m$-compatibility of theta structures as
in \cite[$\S$5.3]{ca05b}.
\begin{lemma}
\label{prodcompat}
Assume that $\Theta_n$ is symmetric.
Then $\Theta_{nm}$ is $m$-compatible with $\Theta_n$.
\end{lemma}
\begin{proof}
Let $\epsilon_n$, $\epsilon_m$,
$E_n$, $E_m$, $\eta_n$, $\eta_m$, $H_n$ and $H_m$ be defined as in \cite[$\S$5.3]{ca05b}.
Note that by the definition of $\Theta_{nm}$ there is an equality
$\Theta_{nm} \circ E_m = \epsilon_m \circ \Theta_n$.
It remains to check that
$\eta_m \circ \Theta_{nm} = \Theta_n \circ H_m$.
In other words, we have to prove that
\begin{eqnarray}
\label{wichtigegl}
\Theta_n \circ H_m \circ E_m = \eta_m \circ \epsilon_m \circ \Theta_n
\quad \mbox{and} \quad
\Theta_n \circ H_m \circ E_n = \eta_m \circ \epsilon_n \circ \Theta_m.
\end{eqnarray}
Using the definition we compute $H_m(E_n(\alpha,x,l))=(\alpha^{nm},0,1)$.
As $(n,m)=1$, it follows that the image of $\eta_m \circ \epsilon_n$
is contained in $\xg_{m,R}$.
Hence the right hand equation in (\ref{wichtigegl}) is a consequence of the
$\xg_{m,R}$-equivariance of $\Theta_n$ and $\Theta_m$.
It remains to prove the left hand equation.
We have
\begin{eqnarray}
\label{eqthree}
\eta_m \circ \epsilon_m =  \delta_m \quad \mbox{and} \quad
H_m \circ E_m  =  D_m
\end{eqnarray}
where $D_m$ denotes the map $G(Z_n) \rightarrow  G(Z_n),(\alpha,x,l) \mapsto
(\alpha^{m^2},mx,l^m)$ and $\delta_m:G(\pol^n) \rightarrow G(\pol^n)$ is given by
\[
g \mapsto g^{(m^2+m)/2} \cdot \delta_{-1}(g)^{(m^2-m)/2}.
\]
Here $\delta_{-1}$ is defined as in Section \ref{thetasym}.
The right hand equation in (\ref{eqthree})
follows by expanding the definitions.
The left hand equation in (\ref{eqthree}) is proven in
\cite[$\S$2, Prop.5]{mu66}.
A straight forward calculation yields that for all $g \in G(Z_n)$ one has
\begin{eqnarray}
\label{eqdef}
D_m(g)=g^{(m^2+m)/2} \cdot D_{-1}(g)^{(m^2-m)/2}.
\end{eqnarray}
The left hand equality in (\ref{wichtigegl})
is implied by the equalities (\ref{eqthree}) and (\ref{eqdef}) using the assumption that
$\Theta_n$ is symmetric,
i.e. equation (\ref{sym}) holds.
This completes the proof of the lemma.
\end{proof}

\subsection{Descent of theta structures by isogeny}

In this section we prove some lemma which forms an important
ingredient of the proof of Theorem \ref{kop}.
The lemma is about special theta relations which are induced by descent along
isogenies.
A proof of this key lemma in terms of algebraic theta
functions is absent from the literature.
In the following we use the notion of
compatibility as defined in \cite[$\S$5.2-5.3]{ca05b}.
\newline\indent
Let $R$ be a local ring, and let 
$\pi_A:A \rightarrow \spec{R}$ and $\pi_B:B \rightarrow \spec{R}$ be abelian
schemes of relative dimension $g$.
We set $Z_n= (\xz / n \xz)^g_R$ for an integer $n \geq 1$.
As usual, we consider $Z_n$ as embedded in $Z_{mn}$ via the morphism that maps
component-wise $1 \mapsto n$.
Let $\emm$ be an ample
symmetric line bundle on $B$.
Suppose that we are given $2$-compatible theta structures $\Sigma_j : G(Z_{jm})
\stackrel{\sim}{\rightarrow} G(\emm^j)$ for some $m \geq 1$, where
$j \in \{ 1,2 \}$.
Let $F:A \rightarrow B$ be an isogeny of degree $d^g$.
Assume that there exists an ample symmetric line bundle
$\pol$ on $A$ such that $F^* \emm \cong
\pol$. Now assume that we are given $2$-compatible theta structures $\Theta_j:
G(Z_{jmd}) \stackrel{\sim}{\rightarrow} G(\pol^j)$
such that $\Theta_{j}$ and $\Sigma_{j}$ are $F$-compatible.
By general theory there exist theta group equivariant isomorphisms
\[
\mu_j:\pi_{A,*} \pol^j \stackrel{\sim}{\rightarrow} V(Z_{jmd})
\quad \mbox{and} \quad
\gamma_j:\pi_{B,*} \emm^j \stackrel{\sim}{\rightarrow} V(Z_{jm}).
\]
Suppose that we have chosen rigidifications of $\pol$ and $\emm$. This defines,
by means of $\mu_j$ and $\gamma_j$, theta functions $q_{ \emm^j } \in V(Z_{jm})$
and $q_{ \pol^{j} } \in V(Z_{jmd})$ (see \cite[$\S$1]{mu66} \cite[$\S$5.1]{ca05b}).
Here we denote the module of algebraic theta functions
by $V(Z_n)=\underline{\mathrm{Hom}}(Z_n, \mathcal{O}_R)$ for an
integer $n \geq 1$.
The following lemma generalizes \cite[Lem.6.4]{ca05b}.
\begin{lemma}
\label{descentone}
There exists a $\lambda \in R^*$ such that for all $x \in Z_{2m}$ one has
\[
q_{ \emm^2 }(x) = \lambda q_{ \pol^{2} }(x).
\]
\end{lemma}
\begin{proof}
By Mumford's \emph{$2$-Multiplication Formula} \cite[$\S$3]{mu66}
there exists a $\lambda \in
R^*$ such that for all $z \in Z_m$ and $x \in Z_{2m}$ we have
\[
(\xone \star \delta_z)(x)= \lambda \sum_{ y \in x + Z_m} \delta_z(x-y)
q_{\emm^2}(y) = \lambda q_{\emm^2}(x-z).
\]
Here $\xone$ denotes the finite theta function which takes the value $1$
on all of $Z_{m}$.
The \emph{Isogeny Theorem} \cite[$\S$1,Th.4]{mu66}
implies that there exists a $\lambda \in R^*$
such that for $x \in Z_{2dm}$ we have
\[
F^*( \xone \star \delta_z)(x)=
\left\{
\begin{array}{c@{, \quad }c}
\lambda q_{\emm^2}(x-z) & x \in Z_{2m} \\
0 & \mathrm{else}
\end{array}
\right.
\]
Also there exists a $\lambda_1,\lambda_2 \in R^*$ such that for $x \in Z_{md}$ we have
\[
F^*( \xone )(x)=
\left\{
\begin{array}{c@{, \quad }c}
\lambda_1 & x \in Z_{m} \\
0 & \mathrm{else}
\end{array}
\right.
\quad \mbox{and} \quad
F^*( \delta_z )= \lambda_2 \delta_z.
\]
Again by Mumford's multiplication formula there exists a $\lambda
\in R^*$ such that for all $x \in Z_{2md}$ we have
\begin{eqnarray*}
& \big( F^*( \xone) \star F^*(\delta_z) \big) (x)
& = \lambda \sum_{ y \in x + Z_{md}} F^*( \xone )(x+y) \delta_z(x-y)
q_{\pol^{2d}}(y) \\
& & = \left\{
\begin{array}{c@{, \quad }c}
\lambda q_{\pol^{2}}(x-z) & x \in Z_{2m} \\
0 & \mathrm{else}
\end{array}
\right.
\end{eqnarray*}
The Lemma now follows from the observation that
$F^*( \xone) \star F^*(\delta_z)$ and $F^*( \xone \star \delta_z)$
differ by a unit.
\end{proof}

\subsection{Products of abelian varieties with theta structure}

In this section we prove the existence of finite products abelian varieties
with theta structures. This kind of product is needed in
the proof of the $3$-multiplication formula.
\newline\indent 
Let $A_1, \ldots, A_n$ be abelian schemes over a ring $R$.
Assume we are given a line bundle $\pol_i$ on $A_i$ and a theta
structure $\Theta_i$ of type $K_i$ for $(A,\pol_i)$ for all $i=1,\ldots,n$.
We set
\[
A= \prod_{i=1}^n A_i,
\quad K= \prod_{i=1}^n K_i
\quad \mbox{and} \quad
\pol= \bigotimes_{i=1}^n p_i^* \pol_i
\]
where $p_i:A \rightarrow A_i$ denotes the projection on the $i$-th
factor.
\begin{lemma}
\label{prodtheta}
There exists a natural product theta structure of type $K$ for $(A, \pol)$
depending on the theta structures $\Theta_i$, where $i=1, \ldots n$.
\end{lemma}
\begin{proof}
We remark that there exists a canonical isomorphism
$\prod_{i=1}^n H(\pol_i) \stackrel{\sim}{\rightarrow} H(\pol)$.
Consider the morphism
$\varphi:\prod_{i=1}^n G(\pol_i) \rightarrow G(\pol)$
given by $(x_i,\psi_i) \mapsto \big( (x_i)_{i=1 \ldots n}, \otimes_{i=1}^n p_i^*
\psi_i \big)$.
Note that
\[
\bigotimes_{i=1}^n p_i^* T_{x_i}^* \pol_i
= \bigotimes_{i=1}^n  T_{(x_1, \ldots, x_n)}^* p_i^* \pol_i
= T_{(x_1, \ldots, x_n)}^* \pol.
\]
Obviously, $C=\mathrm{ker}(\varphi)
\subseteq \xg_{m,R}^n$.
We claim that $\varphi$ is surjective.
Consider the diagram
\begin{eqnarray*}
\xymatrix{ \prod_{i=1}^n G(\pol_i)  \ar@{->}[rr]^{\varphi}
  \ar@{->}[d]_{\pi_1 \times \ldots \times \pi_n} & &
  G(\pol) \ar@{->}[d]^{\pi} \\
  \prod_{i=1}^n H(\pol_i)
  \ar@{->}[rr]^{\mathrm{can}} \ar@/_1pc/[u]_s & &
  H(\pol) },
\end{eqnarray*}
where $\pi_i$ ($i=1, \ldots, n$) and $\pi$ denote the natural
projections. Let $s$ be the canonical section of $\pi_1 \times \ldots \times \pi_n$
induced by the theta structures $\Theta_i$.
We have $\pi \circ \varphi = \pi_1 \times \ldots
\times \pi_n$ (up to canonical isomorphism).
As a consequence we have $\pi \circ \varphi \circ s
= (\pi_1 \times \ldots \times \pi_n) \circ s = \mathrm{id}$.
We conclude that
$\pi \circ \varphi \circ s \circ \pi = \pi$.
Let $g \in G(\pol)$. Then by the latter equality the group
element $\varphi(s(\pi(g)))$ differs from $g$ by a
unit. Hence the morphism $\varphi$ maps a suitable multiple
of $s( \pi(g) )$ to $g$. This implies the
surjectivity of $\varphi$ and proves our claim.
\newline\indent
Analogously, one defines a surjective morphism
$\Phi : \prod_{i=1}^n G(K_i) \rightarrow G(K)$
having kernel equal to $C$.
Let $\tilde{\varphi}$ and $\tilde{\Phi}$ denote the induced isomorphisms
\[
\left( \prod_{i=1}^n G(\pol_i) \right) / C \stackrel{\sim}{\rightarrow} G(\pol)
\quad \mbox{and} \quad
\left( \prod_{i=1}^n G(K_i) \right) / C \stackrel{\sim}{\rightarrow} G(K).
\]
The theta structure $\Theta$ whose existence is claimed in the lemma
is given by the composed isomorphism
\[
G(K) \stackrel{\tilde{\Phi}^{-1}}{\longrightarrow}
\left( \prod_{i=1}^n G(K_i) \right) /C
\quad
\stackrel{\Theta_1 \times \ldots \times \Theta_n}{\longrightarrow}
\quad
\left( \prod_{i=1}^n G(\pol_i) \right) /C
\stackrel{\tilde{\varphi}}{\longrightarrow} G(\pol).
\] 
This completes the proof of the lemma.
\end{proof}

\subsection{An algebraic proof of the $3$-multiplication formula}
\label{mult3}

In the following we give
a \emph{$3$-multiplication formula} for algebraic theta functions
in the context of Mumford's theory \cite{mu66}.
Our method of proof extends to an arbitrary $n$-product of
algebraic theta functions.
For this reason it seems to be instructive
to give a detailed proof in terms of Mumford's algebraic theta
functions.
The following proof generalizes in a straight forward manner
Mumford's proof of his \emph{$2$-multiplication formula} \cite[$\S3$]{mu66}.
The classical complex analytic $3$-multiplication formula does not
apply in our case because we are working in an arithmetic setting.
The theory of algebraic theta functions allows us to keep track of the
reduction modulo the prime $3$. Let us remind the reader, that
our aim is to use the $3$-multiplication
formula in order to lift theta null points from the special fiber
to characteristic $0$.
For the proof of the complex analytic $3$-multiplication formula
we refer to \cite[Ch.7.6]{lb03}.
\newline\indent
Let $A$ be an abelian scheme over a local ring $R$ and let
$\xi$ denote the isogeny $A^3 \rightarrow A^3$ given by
\[
(x_1,x_2,x_3) \mapsto (x_1-2x_2,x_1+x_2-x_3,x_1+x_2+x_3)
\]
Assume we are given an ample line bundle $\pol$ on $A$
and theta structures $\Theta_i$ of type $K_i$ for $\pol^i$
where $i \in I= \{ 1,2,3,6 \}$.
We assume that the theta structures $\Theta_i$, $i \in I$,
are compatible in the sense of \cite[$\S$5.3]{ca05b}.
We set
\[
\emm_{i,j,l}=p_1^* \pol^i \otimes p_2^* \pol^j \otimes
p_3^* \pol^l,
\]
where $p_r:A^3 \rightarrow A$, $r=1,2,3$, is the projection on the $r$-th factor,
and $K_{i,j,l}=K_i \times K_j \times K_l$
for $i,j,l \in I$.
By Lemma \ref{prodtheta}
there exist product theta
structures $\Theta_{1,1,1}$ and $\Theta_{3,6,2}$
of type $K_{1,1,1}$ and $K_{3,6,2}$ for $\emm_{1,1,1}$ and $\emm_{3,6,2}$,
respectively,
depending on the theta structures $\Theta_i$ where $i \in I$.
\begin{proposition}
\label{pullbackxi}
There exists an isomorphism
\begin{eqnarray}
\label{defiso}
\xi^* \emm_{1,1,1} \stackrel{\sim}{\rightarrow} \emm_{3,6,2}.
\end{eqnarray}
\end{proposition}
\begin{proof}
Let $b=(b_1,b_2) \in A^2$ and $a \in A$. We define
\[
s_1:A^2 \rightarrow A^3,(x_1,x_2) \mapsto (a,x_1,x_2)
\quad
\mbox{and} \quad
s_2:A \rightarrow A^3,x \mapsto (x,b_1,b_2).
\]
One computes
\begin{eqnarray*}
s_2^* \emm_{3,6,2} = s_2^* p_1^* \pol^3 \otimes s_2^* p_2^* \pol^6
\otimes s_2^* p_3^* \pol^2
= (p_1 \circ s_2)^* \pol^3 \otimes (p_2 \circ s_2)^* \pol^6 \otimes
(p_3 \circ s_2)^* \pol^2 = \pol^3.
\end{eqnarray*}
and
\begin{eqnarray*}
s_2^* \xi^* \emm_{1,1,1} & = & (p_1 \circ \xi \circ s_2)^* \pol
\otimes (p_2 \circ \xi \circ s_2)^* \pol \otimes (p_3 \circ \xi \circ
s_2)^* \pol \\
& = & T_{-2b_1}^* \pol \otimes T_{b_1-b_2}^* \pol
\otimes T_{b_1+b_2}^* \pol = \pol^3.
\end{eqnarray*}
The latter equality follows by the Theorem of the Square.
Now take $a=0_A$ where $0_A$ denotes the zero section of $A$.
Let $p_{23}:A^3 \rightarrow A^2$ be the projection on the $2$-nd and
$3$-rd factor and let $\tilde{p}_m:A^2 \rightarrow A$ denote the projection on the $m$-th
factor ($m=1,2$).
We have
\begin{eqnarray*}
s_1^* \emm_{3,6,2}
& = & ( p_1 \circ s_1)^* \pol^3 \otimes s_1^*( p_2^* \pol^6
\otimes p_3^* \pol^2) \\
& = & ( p_1 \circ s_1)^* \pol^3
\otimes (p_{23} \circ s_1)^* (\tilde{p}_1^* \pol^6
\otimes \tilde{p}_2^* \pol^2)
= \tilde{p}_1^* \pol^6
\otimes \tilde{p}_2^* \pol^2.
\end{eqnarray*}
By \cite[$\S$3,Prop.1]{mu66} we conclude that
\begin{eqnarray*}
s_1^* \xi^* \emm_{1,1,1} &
= & (p_1 \circ \xi \circ s_1)^* \pol
\otimes (p_{23} \circ \xi \circ s_1)^* (\tilde{p}_1^* \pol
\otimes \tilde{p}_2^* \pol) \\
& & = \tilde{p}_1^* [2]^* [-1]^* \pol
\otimes (\tilde{p}_1^* \pol
\otimes \tilde{p}_2^* \pol)^2 = \tilde{p}_1^* \pol^{6}
\otimes \tilde{p}_2^* \pol^{2}.
\end{eqnarray*}
The latter equality follows by the symmetry of $\pol$.
Note that $p_{23} \circ \xi \circ s_1$ equals the isogeny
used in \cite[$\S$3,Prop.1]{mu66}.
The proposition now follows by applying the Seesaw Principle.
\end{proof}
\begin{lemma}
\label{prepmult}
The theta structure $\Theta_{3,6,2}$ is $\xi$-compatible with
$\Theta_{1,1,1}$.
\end{lemma}
\begin{proof}
We have to check the compatibility assumptions of
\cite[$\S$5.2]{ca05a}.
We have already shown in Proposition \ref{pullbackxi} that
there exists an isomorphism
$\alpha:\xi^* \emm_{1,1,1} \stackrel{\sim}{\rightarrow} \emm_{3,6,2}$.
Let $\tau$ be the morphism $A \rightarrow A^3,x \mapsto (2x,x,3x)$.
The kernel of $\xi$ is given by the restriction of $\tau$ to $A[6]$.
In the following we will identify the groups $K_6 \times K_6^D$,
$K_{1,1,1} \times K_{1,1,1}^D$
and $K_{3,6,2} \times K_{3,6,2}^D$ with their
images under the Lagrangian decompositions induced by
the theta structures $\Theta_6$, $\Theta_{1,1,1}$ and
$\Theta_{3,6,2}$, respectively.
Note that $A[6]$ is contained in the image of $K_6 \times K_6^D$ under
$\tau$.
By the compatibility assumptions we have
$\tau \big( K_6 ) \subseteq K_{3,6,2}$ and
$\tau \big( K_6^D \big) \subseteq K_{3,6,2}^D$.
We conclude that
condition ($\dagger$) of \cite[$\S$5.2]{ca05a} is satisfied with
$Z_1= \tau \big( A[6] \cap K_6 \big)$ and
$Z_2= \tau \big( A[6] \cap K_6^D \big)$.
The isomorphism (\ref{defiso}) gives rise to a subgroup
$\tilde{K} \leq G(\emm_{3,6,2})$
lifting the kernel of $\xi$. Let $G(\emm_{3,6,2})^*$ denote the centralizer
of $\tilde{K}$ in $G(\emm_{3,6,2})$.
By \cite[$\S$1,Prop.2]{mu66} we have
\begin{eqnarray}
\label{mumeq}
G(\emm_{3,6,2})^* = \left\{ g \in G(\emm_{3,6,2})| \xi \big( \pi_{3,6,2}(g)
  \big) \in A^3[\emm_{1,1,1}] \right\}
\end{eqnarray}
where $\pi_{3,6,2}:G(\emm_{3,6,2})
\rightarrow A^3[\emm_{3,6,2}]$ is the natural projection.
Here we denote
\[
A^3[\emm_{i,j,l}]= \{ x \in A^3 \vert T_x^*
\emm_{i,j,l} \cong \emm_{i,j,l} \}
\]
for all $i,j,l \in I$.
Because of the equality (\ref{mumeq}) we have
\[
Z_1^{\bot}=\{ (x,y,z) \in K_{3,6,2} | \xi(x,y,z) \in K_{1,1,1} \}
\quad \mbox{and} \quad
Z_2^{\bot}=\{ (x,y,z) \in K_{3,6,2}^D | \xi(x,y,z) \in K_{1,1,1}^D \}
\]
(notation as in \cite[$\S$5.2]{ca05a}).
Obviously the isogeny $\xi$ induces a surjective morphism
$\sigma: Z_1^{\bot} \rightarrow K_{1,1,1}$ having kernel $Z_1$.
Let $\sigma_1$ be the inverse of the isomorphism $Z_1^{\bot} / Z_1
\stackrel{\sim}{\rightarrow} K_{1,1,1}$ induced by $\sigma$.
Let $\sigma_2$ be defined as in \cite[$\S$5.2]{ca05b}.
It remains to check the commutativity of the following diagram
\begin{eqnarray}
\label{commdia}
\xymatrix{
G( \emm_{3,6,2} )^* / \tilde{K} \ar@{<-}[rr]^{\Theta_{3,6,2}}
\ar@{->}[d]^{\mathrm{can},\xi} & & \xg_{m,R} \times Z_1^{\bot} / Z_1
\times Z_2^{\bot} / Z_2 \ar@{<-}[d]^{\mathrm{id} \times \sigma_1
\times \sigma_2} \\
G(\emm_{1,1,1}) \ar@{<-}[rr]^{\Theta_{1,1,1}} & & \xg_{m,R} \times
K_{1,1,1} \times K_{1,1,1}^D}
\end{eqnarray}
where the left hand vertical morphism is defined as in the proof of
\cite[$\S$1,Prop.2]{mu66}.
We claim that the group $ \xg_{m,R} \times Z_1^{\bot} / Z_1
\times Z_2^{\bot} / Z_2$ is generated by $\xg_{m,R}$ and elements of the form
\[
(1,2x,x,3x,l^2,l,l^3), \quad (1,2x,x,-3x,l^2,l,l^{-3}) \quad
\mbox{and} \quad (1,2x,-2x,0,l^2,l^{-2},1)
\]
where $(x,l) \in K_6 \times K_6^D$.
Let $\xi'$ denote the isogeny $A^3 \rightarrow A^3$ given by
\[
(x_1,x_2,x_3) \mapsto (2x_1+2x_2+2x_3,-2x_1+x_2+x_3,-3x_2+3x_3).
\]
Assume we are given an element $(1,x,l)$ of
$\xg_{m,R} \times Z_1^{\bot}  \times Z_2^{\bot} $.
We denote
$\xi (x,l) =(\bar{x},\bar{l})$.
Choose $\tilde{x} \in K_{6,6,6}$ and $\tilde{l} \in K_{6,6,6}^D$
such that $[6](\tilde{x},\tilde{l})=(\bar{x},\bar{l})$.
One verifies that $\xi \circ \xi' = [6]$ and hence the element
$\xi' (\tilde{x},\tilde{l}) \in K_{3,6,2} \times K_{3,6,2}^D$
differs from $(x,l)$ by an element of $Z_1 \times Z_2$.
This implies the above claim.
\newline\indent
In the sequel we will prove the commutativity of the diagram (\ref{commdia}) for
elements of the form $(1,2x,x,3x,l^2,l,l^3)$.
The proof for elements of the form
\[
(1,2x,x,-3x,l^2,l,l^{-3}) \quad \mbox{and} \quad
(1,2x,-2x,0,l^2,l^{-2},1)
\]
goes analogously and is left to the reader.
We define
\[
\iota:G(K_6) \rightarrow G(K_{3,6,2}), ( \alpha, x,l)
\mapsto ( \alpha^6, 2x,x,3x , l^2, l, l^3)
\]
and set
$\kappa=\Theta_{3,6,2} \circ \iota \circ \Theta_6^{-1}$.
Let
$G(K_6)^{\sharp}=\iota^{-1} G(K_{3,6,2})^*$ and
$G(\pol^6)^{\sharp}=\Theta_6 \big( G(K_6)^{\sharp} \big)$.
We define
$\varphi_3:G(\pol) \rightarrow G(\emm_{1,1,1})$ and
$\Phi_3:G(K_1) \rightarrow G(K_{1,1,1})$ to
be the restriction on the $3$-rd factor of the morphism $\varphi$ and $\Phi$ introduced
in the proof of Lemma \ref{prodtheta}.  
It is readily checked that the following diagram (dotted arrows
ignored) is commutative
\begin{eqnarray*}
\xymatrix{
G(\pol^6)^{\sharp} \ar@{->}[dd]_{\eta_6} \ar@{<-}[rr]^{\Theta_6}
\ar@{->}[rd]^{\kappa} & &
G(K_6)^{\sharp} \ar@{->}[dd]_(0.3){H_6} \ar@{->}[rd]^{\iota} & \\
& G(\emm_{3,6,2})^* / \tilde{K} \ar@{<-}[rr]_(0.6){\Theta_{3,6,2}}
\ar@{.>}[dd]^(0.3){\mathrm{can},\xi}
& & \xg_{m,R} \times Z_1^{\bot} / Z_1 \times Z_2^{\bot} / Z_2
\ar@{<-}[dd]^{\mathrm{id} \times \sigma_1 \times \sigma_2} \\
G(\pol) \ar@{<-}[rr]^(0.3){\Theta_1} \ar@{->}[dr]^{\varphi_3} & &
G(K_1) \ar@{->}[rd]^{\Phi_3} & \\
& G(\emm_{1,1,1}) \ar@{<-}[rr]^{\Theta_{1,1,1}} & &
\xg_{m,R} \times K_{1,1,1} \times K_{1,1,1}^D.
}
\end{eqnarray*}
Here $\eta_6$ and $H_6$ are defined as in \cite[$\S$5.3]{ca05b}. Note
that the upper left square is commutative since $\Theta_6$ and
$\Theta_1$ are assumed to be $6$-compatible.
\newline\indent
In order to show that the diagram (\ref{commdia}) is commutative
on the subset of elements of the form $(1,2x,x,3x,l^2,l,l^3)$
it suffices to prove that the following diagram commutes
\begin{eqnarray*}
\xymatrix{
G(\pol^6)^{\sharp} \ar@{->}[rr]^{\kappa} \ar@{->}[d]_{\eta_6} & &
G(\emm_{3,6,2})^* / \tilde{K}  \ar@{->}[d]^{\mathrm{can},\xi} \\
G(\pol) \ar@{->}[rr]^{\varphi_3} & & G(\emm_{1,1,1}). 
}
\end{eqnarray*}
Consider the commutative diagram
\begin{eqnarray*}
\xymatrix{
A^3 \ar@{->}[r]^{\xi} \ar@{<-}[d]_{\tau} & A^3 \ar@{<-}[d]^{\mathrm{i}_3} \\
A \ar@{->}[r]^{[6]} & A }
\end{eqnarray*}
where $\mathrm{i}_3(x)= (0,0,x)$.
There exist isomorphisms
\begin{eqnarray*}
\beta:[6]^* \pol \stackrel{\sim}{\rightarrow} \pol^{36}
\quad \mbox{and} \quad
\gamma:\mathrm{i}_3^* \emm_{1,1,1} \stackrel{\sim}{\rightarrow} \pol.
\end{eqnarray*}
The existence of the isomorphism $\beta$ is implied by the symmetry of $\pol$.
Consider the isomorphism $\delta$ given by the composition
\begin{eqnarray*}
\tau^* \emm_{3,6,2} \stackrel{\tau^* \alpha^{-1}}{\longrightarrow}
\tau^* \xi^* \emm_{1,1,1} = (\xi \circ \tau)^* \emm_{1,1,1}
= (i_3 \circ [6])^* \emm_{1,1,1}
= [6]^* i_3^* \emm_{1,1,1} \stackrel{[6]^* \gamma}{\longrightarrow}
[6]^* \pol \stackrel{\beta}{\rightarrow} \pol^{36},
\end{eqnarray*}
where $\alpha$ is as above.
The isomorphism $\delta$ induces a morphism
\[
\mathrm{can},\tau:G(\emm_{3,6,2}) \rightarrow G(\pol^{36}),
( x , \psi ) \mapsto \big( p_2(x) , T_{p_2(x)}^* \delta \circ \tau^* \psi
\circ \delta^{-1} \big)
\]
where $p_2:A^3 \rightarrow A$ denotes the projection on the second factor.
We claim that the following diagram is commutative
\begin{eqnarray}
\label{epsilondia}
\xymatrix{
G(\pol^6) \ar@{->}[r]^{\kappa} \ar@{->}[d]_{\epsilon_6}
& G(\emm_{3,6,2}) \ar@{->}[dl]^{\mathrm{can},\tau} \\
G(\pol^{36}) & }
\end{eqnarray}
where $\epsilon_6$ is defined as in \cite[$\S$5.3]{ca05b}.
Let $g=\big( (x,l),\psi \big) \in G(\pol^6)$ and
$h=\Theta_6^{-1}(g)$. By definition we have
\[
\iota ( h ) = \Phi \big( H_2(h),h,H_3(h) \big)
\]
and hence
\begin{eqnarray*}
& \kappa(g)= \varphi \big( \eta_2(g), g , \eta_3(g) \big)=
\big( (2x,x,3x,l^2,l,l^3) , p_1^* \eta_2( \psi) \otimes p_2^* \psi
\otimes p_3^* \eta_3 (\psi) \big).&
\end{eqnarray*}
The image of $\kappa(g)$ under the canonical morphism induced by
$\delta$ is given by
\[
\Big( x, T_x^* \delta \circ \tau^* \big(
 p_1^* \eta_2( \psi) \otimes p_2^* \psi
\otimes p_3^* \eta_3 (\psi) \big) \circ \delta^{-1} \Big).
\]
Choose isomorphisms
$\rho_2:[2]^* \pol \stackrel{\sim}{\rightarrow} \pol^4$
and $\rho_3: [3]^* \pol \stackrel{\sim}{\rightarrow} \pol^9$.
Consider the composed isomorphism $\delta'$ given by
\begin{eqnarray*}
& \tau^* \emm_{3,6,2}
= \tau^* ( p_1^* \pol^3 \otimes p_2^* \pol^6 \otimes p_3^* \pol^2 )
= (p_1 \circ \tau)^* \pol^3 \otimes (p_2 \circ \tau)^* \pol^6
\otimes (p_3 \circ \tau)^* \pol^2 & \\
& = [2]^* \pol^3 \otimes \pol^6 \otimes [3]^* \pol^2
\stackrel{\rho_2 \otimes \mathrm{id} \otimes \rho_3}{\longrightarrow}
\pol^{36}.& 
\end{eqnarray*}
The isomorphism $\delta'$ differs from $\delta$ by a unit.
Thus we have
\begin{eqnarray*}
\lefteqn{T_x^* \delta \circ \tau^* \big(
 p_1^* \eta_2( \psi) \otimes p_2^* \psi
\otimes p_3^* \eta_3 (\psi) \big) \circ \delta^{-1}} \\
&&= T_x^* \delta' \circ \tau^* \big(
 p_1^* \eta_2( \psi) \otimes p_2^* \psi
\otimes p_3^* \eta_3 (\psi) \big) \circ (\delta')^{-1} \\
& & = T_x^* \delta' \circ \big(
(p_1 \circ \tau)^* \eta_2( \psi) \otimes (p_2 \circ \tau)^* \psi
\otimes (p_3 \circ \tau)^* \eta_3 (\psi) \big)
\circ (\delta')^{-1} \\
& & = T_x^* \delta' \circ \big(
[2]^* \eta_2( \psi) \otimes \psi
\otimes [3]^* \eta_3 (\psi) \big)
\circ (\delta')^{-1} \\
& & = \big( T_x^* \rho_2 \circ [2]^* \eta_2( \psi) \circ \rho_2^{-1}
\big) \otimes \psi \otimes \big( T_x^* \rho_3 \circ  [3]^* \eta_3
(\psi) \circ \rho_3^{-1} \big) \\
&&= \epsilon_2(\psi) \otimes \psi \otimes \epsilon_3(\psi)
= \epsilon_6(\psi).
\end{eqnarray*}
This proves our claim, i.e. the commutativity of diagram
(\ref{epsilondia}).
The isomorphism $\gamma$ induces a morphism
\[
\mathrm{can},i_3:G(\emm_{1,1,1}) \rightarrow G(\pol),
( x , \psi ) \mapsto \big( p_3(x) , T_{p_3(x)}^* \gamma \circ i_3^* \psi
\circ \gamma^{-1} \big)
\]
where $p_3:A^3 \rightarrow A$ denotes the projection on the $3$rd factor.
Consider the diagram
\begin{eqnarray*}
\xymatrix{
G(\emm_{3,6,2})^* \ar@{->}[rr]
\ar@{->}[dd]^{\mathrm{can},\tau} \ar@{<-}[rrd]^{\kappa} &
& G(\emm_{3,6,2})^* / \tilde{K} \ar@{->}[rr]^{\mathrm{can},\xi} &
& G(\emm_{1,1,1}) \ar@/^/[dd]^{\mathrm{can},i_3} \\
& & G(\pol^6)^{\sharp} \ar@{->}[dll]_{\epsilon_6}
\ar@{->}[drr]^{\eta_6} &
& & \\
G(\pol^{36}) \ar@{<-^)}[r] &
G(\pol^{36})^* \ar@{->}[r] &
G(\pol^{36})^*/ \widetilde{A[6]} \ar@{->}[rr]_{\mathrm{can},[6]} &
& G(\pol) \ar@/^/[uu]^{\varphi_3}.
}
\end{eqnarray*}
Here $G(\pol^{36})^*$ denotes the centralizer of the lifted
subgroup $\widetilde{A[6]}$ in $G(\pol^{36})$.
By the above discussion the left hand triangle is commutative.
By the same reasoning as above it follows that the composed morphism
\[
G(\pol) \stackrel{\varphi_3}{\longrightarrow} G(\emm_{1,1,1})
\stackrel{\mathrm{can},i_3}{\longrightarrow} G(\pol)
\]
equals the identity.
This implies that the canonical morphism induced by $\gamma$ is
surjective.
As a consequence the commutativity of diagram (\ref{commdia})
is equivalent to the commutativity of the following diagram
\begin{eqnarray}
\label{eqdia}
\xymatrix{
G(\pol^6)^{\sharp} \ar@{->}[rr]^{\kappa} \ar@{->}[d]^{\eta_6} & &
G(\emm_{3,6,2})^* / \tilde{K} \ar@{->}[d]^{\mathrm{can},\xi}  \\
G( \pol ) \ar@{<-}[rr]^{\mathrm{can},i_3} & & G(\emm_{1,1,1}).
}
\end{eqnarray}
Let
\[
g \in G(\pol^6)^{\sharp} \quad \mbox{and} \quad
\kappa(g)=(x, \psi).
\]
By definition the
image of $\kappa(g)$ under the canonical morphism induced by $\delta$
is given by
\begin{eqnarray*}
\Big( p_2(x), T_{p_2(x)}^* \delta \circ \tau^* \psi \circ \delta^{-1} \Big).
\end{eqnarray*}
Since $\tau \big( p_2(x) \big)=x$ it follows that
\begin{eqnarray*}
\lefteqn{  T_{p_2(x)}^* \delta \circ \tau^* \psi \circ \delta^{-1}
= T_{p_2(x)}^* \big( \beta \circ [6]^* \gamma \circ \tau^*
\alpha^{-1} \big) \circ \tau^* \psi \circ \big( \beta \circ [6]^* \gamma \circ \tau^*
\alpha^{-1} \big)^{-1}} \\
& & =  T_{p_2(x)}^*  \beta \circ T_{p_2(x)}^* [6]^* \gamma \circ
T_{p_2(x)}^* \tau^* \alpha^{-1} \circ \tau^* \psi \circ
\tau^* \alpha \circ [6]^* \gamma^{-1} \circ \beta^{-1} \\
& & =  T_{p_2(x)}^*  \beta \circ [6]^* T_{p_2(6x)}^* \gamma \circ
\tau^* T_{x}^* \alpha^{-1} \circ \tau^* \psi \circ
\tau^* \alpha \circ [6]^* \gamma^{-1} \circ \beta^{-1} \\
& & = T_{p_2(x)}^*  \beta \circ [6]^* T_{p_2(6x)}^* \gamma \circ
\tau^* \big( T_{x}^* \alpha^{-1} \circ \psi \circ
\alpha \big) \circ [6]^* \gamma^{-1} \circ \beta^{-1} \\
& & = T_{p_2(x)}^*  \beta \circ [6]^* T_{p_2(6x)}^* \gamma \circ
\tau^* \xi^* \psi' \circ [6]^* \gamma^{-1} \circ \beta^{-1}
\end{eqnarray*}
where $\xi^* \psi'= T_{x}^* \alpha^{-1} \circ \psi \circ
\alpha$.
Note that such an isomorphism $\psi'$ exists since $\kappa(g) \in G(\emm_{3,6,2})^*$.
We remark that the pair
$\big( \xi (x) , \psi' \big)
\in G(\emm_{1,1,1})$
is the image of $\kappa(g)$ under the canonical morphism
induced by $\alpha$.
Continuing the above calculation we get
\begin{eqnarray*}
\lefteqn{  T_{p_2(x)}^* \delta \circ \tau^* \psi \circ \delta^{-1}
= T_{p_2(x)}^*  \beta \circ [6]^* T_{p_2(6x)}^* \gamma \circ
\tau^* \xi^* \psi' \circ [6]^* \gamma^{-1} \circ \beta^{-1} } \\
& & = T_{p_2(x)}^*  \beta \circ [6]^* T_{p_2(6x)}^* \gamma \circ
[6]^* i_3^* \psi' \circ [6]^* \gamma^{-1} \circ \beta^{-1} \\
& & = T_{p_2(x)}^*  \beta \circ [6]^* \big( T_{p_2(6x)}^* \gamma \circ
i_3^* \psi' \circ \gamma^{-1} \big) \circ \beta^{-1} \\
& & = T_{p_2(x)}^*  \beta \circ [6]^* \big( T_{p_3(\xi(x))}^* \gamma \circ
i_3^* \psi' \circ \gamma^{-1} \big) \circ \beta^{-1}.
\end{eqnarray*}
By definition the pair
\[
g'=\left( p_3 \big( \xi(x) \big) , T_{p_3(\xi(x))}^* \gamma \circ
i_3^* \psi' \circ \gamma^{-1} \right)
\]
is the image of $\big( \xi (x) , \psi' \big)$ under the canonical
morphism induced by $\gamma$.
We conclude by the above equality, the commutativity of diagram
(\ref{epsilondia}) and the definition of
$\eta_6$ that $g'= \eta_6(g)$.
Thus we have shown that diagram (\ref{eqdia}) is commutative. As a
consequence diagram (\ref{commdia}) is commutative.
This finishes the proof of the lemma.
\end{proof}
Assume that we have chosen $G(K_{1,1,1})$- and
$G(K_{3,6,2})$-equivariant isomorphisms
\[
\mu_{1,1,1}: \pi_{3,*} \emm_{1,1,1} \stackrel{\sim}{\rightarrow} V(K_{1,1,1})
\quad \mbox{and} \quad
\mu_{3,6,2}: \pi_{3,*} \emm_{3,6,2} \stackrel{\sim}{\rightarrow} V(K_{3,6,2})
\]
where $\pi_3$ denotes the structure morphism of $A^3$.
The following lemma is a generalization of the \emph{Addition Formula}
which is stated in \cite[$\S$3]{mu66}.
We use the intuitively simplified notation introduced in the proof of
Lemma \ref{prepmult}.
\begin{corollary}
\label{addform}
There exists a $\lambda \in R^*$ such that for all $g \in
V(K_{1,1,1})$ we have
\[
\xi^*g(x,y,z)= \left\{
\begin{array}{l@{, \quad}l}
\lambda g(\xi(x,y,z)) & \xi(x,y,z) \in K_{1,1,1} \\
0 & \mathrm{else}
\end{array} \right.
\]
where $(x,y,z) \in K_{3,6,2}$.
\end{corollary}
\begin{proof}
By Lemma \ref{prepmult} we can apply the \emph{Isogeny Theorem} (see
\cite[$\S$1,Th.4]{mu66} \cite[$\S$5.2,Th.5.4]{ca05b}) in order
to obtain the formula given in the lemma.
\end{proof}
Assume that we have chosen $G(K_i)$-equivariant isomorphisms
\[
\mu_i: \pi_* \pol^i \stackrel{\sim}{\rightarrow} G(K_i), \quad i \in I,
\]
where $\pi$ denotes the structure morphism of $A$, and
that we have rigidified the line bundle $\pol$.
This defines theta functions $q_{\pol^i} \in V(K_i)$
(see \cite[$\S$1]{mu66} and \cite[$\S$5.1]{ca05b}).
\newline\indent
Let $\Delta:A
\rightarrow A^3$ the diagonal morphism.
There exists a canonical isomorphism
$\beta:\Delta^* \emm_{1,1,1} \stackrel{\sim}{\rightarrow} \pol^3$.
The following theorem describes the morphism of
$\mathcal{O}_R$-modules
$\varphi$ defined as the composition
\begin{eqnarray*}
\pi_* \pol \otimes \pi_* \pol \otimes \pi_* \pol
\stackrel{\mathrm{can}}{\rightarrow} \pi_{3,*} \emm_{1,1,1}
\stackrel{\mathrm{can}}{\rightarrow} \pi_{3,*} \Delta_*
\Delta^* \emm_{1,1,1} = \pi_* \Delta^* \emm_{1,1,1}
\stackrel{\pi_* \beta}{\longrightarrow} \pi_* \pol^3,
\end{eqnarray*}
where the left hand morphism is the K\"{u}nneth morphism,
in terms of finite theta functions.
\begin{definition}
For $s_1,s_2,s_3 \in \pi_* \pol$ and
$f_1,f_2,f_3 \in V(K_1)$ such that $\mu(s_i)=f_i$ ($i=1,2,3$) we set
\[
f_1 \star f_2 \star f_3 = (\mu_1 \otimes \mu_1 \otimes \mu_1) \big(s_1
\otimes s_2 \otimes s_3 \big).
\]
\end{definition}
We define for $x \in K_3$
\[
G_x = \{ (y,z) \in K_{6,2} \hspace{0.2cm} | \hspace{0.2cm} \xi(x,y,z) \in K_{1,1,1} \}.
\]
\begin{theorem}[$3$-multiplication formula]
\label{threemult}
There exists a $\lambda \in R^*$ such that for all $x \in K_3$ and
$f_1,f_2,f_3 \in V(K_1)$ we have
\[
(f_1 \star f_2 \star f_3)(x)= \lambda  \sum_{ (y,z) \in G_x} f_1(x-2y) f_2(x+y-z)
f_3(x+y+z) q_{ \pol^6}(y) q_{ \pol^2 }(z).
\]
\end{theorem}
\begin{proof}
Consider the commutative diagram
\begin{eqnarray}
\label{didia}
\xymatrix{
A \ar@{->}[d]^{i_1} \ar@{->}[rd]^{\Delta} & \\
A^3 \ar@{->}[r]^{\xi} & A^3
}
\end{eqnarray}
where $i_1:A \rightarrow A^3$ is defined by $x \mapsto (x,0,0)$ and
$\Delta$ is the diagonal morphism.
Note that there exists an isomorphism $\gamma:i_1^* \emm_{3,6,2} \stackrel{\sim}
{\rightarrow} \pol^3$.
By Proposition \ref{pullbackxi} there exists an isomorphism
$\alpha: \xi^* \emm_{1,1,1} \stackrel{\sim}{\rightarrow} \emm_{3,6,2}$.
Because of the commutativity of the diagram (\ref{didia})
the morphism $\varphi$ defined above equals up to a unit the composed morphism
\begin{eqnarray*}
\lefteqn{\pi_* \pol \otimes \pi_* \pol \otimes \pi_* \pol
\stackrel{\mathrm{can}}{\longrightarrow} \pi_{3,*} \emm_{1,1,1}
\stackrel{\mathrm{can}}{\longrightarrow} \pi_{3,*} \Delta_*
\Delta^* \emm_{1,1,1}} \\
& &  = \pi_* \Delta^* \emm_{1,1,1} = \pi_* i_1^* \xi^* \emm_{1,1,1}
\stackrel{\pi_* i_1^* \alpha }{\longrightarrow} \pi_* i_1^*
\emm_{3,6,2}
\stackrel{\pi_* \gamma}{\longrightarrow} \pi_* \pol^3.
\end{eqnarray*}
Passing over from sections to finite theta functions we get a diagram
\begin{eqnarray}
\label{keymorph}
V(K_1) \otimes V(K_1) \otimes V(K_1) \stackrel{\mathrm{can}}{\rightarrow} V(K_{1,1,1})
\stackrel{\xi^*}{\rightarrow} V(K_{3,6,2})
\stackrel{\mathrm{eval}}{\rightarrow} V(K_3).
\end{eqnarray}
The left hand map is defined to be the canonical isomorphism mapping
\[
f_1 \star f_2 \star f_3 \mapsto \tilde{f}_1 \tilde{f}_2
\tilde{f}_3
\]
where $\tilde{f}_i$ is the function on $K_{1,1,1}$ defined by
\[
\tilde{f}_i(x_1,x_2,x_3) = f_i(x_i), \quad i=1,2,3.
\]
The map $\xi^*$ is given by Corollary \ref{addform}.
The right hand $\mathrm{eval}$-map in diagram (\ref{keymorph}) corresponds to the map on sections
which maps a section $s_1 \otimes s_2 \otimes s_3 \in \pi_* \emm_{3,6,2}$
to the section $(s_2)_0 (s_3)_0 s_1$ where $(\cdot)_0$ indicates the
evaluation at zero by means of the chosen rigidification.
The claim now follows by expressing $(s_2)_0$ and $(s_3)_0$ in terms
of theta null values (see \cite[$\S$1,Cor.3]{mu66}).
\end{proof}

\section{Explicit CM construction in characteristic $3$}
\label{algo}

In this section we apply Corollary~\ref{genus2} to the explicit CM 
construction of invariants of ordinary abelian surfaces by canonical lifting 
from characteristic~$3$. The CM algorithm has two main phases:
\begin{itemize}
\item
first (see Section \ref{lift_g2}), the multivariate Newton lifting of a 
given canonical theta null point based on the equations of
Corollary~\ref{genus2} by means of the algorithm of Lercier and 
Lubicz~\cite[Th.2]{ll05},
\item
second (see Section \ref{LLL}), the LLL reconstruction of the defining 
polynomials over $\ZZ$ for the ideal of relations between the canonically 
lifted moduli, following Gaudry et al.~\cite{Gaudryandall}.
\end{itemize}
The existence of the lifting algorithm is a consequence of the following facts.
The ordinary locus at $3$ of the moduli space of abelian varieties with symmetric 
$4$-theta structure, which is constructed in \cite{mu67a}, is smooth. 
The space of pairs of ordinary abelian varieties with symmetric $4$-theta 
structure admitting a compatible isogeny of degree $3^g$, where $g$ is the 
dimension, forms an \'{e}tale covering of the latter space.
\newline\indent
The lifting algorithm applies to a rationally parametrized moduli space $X$ 
over $\ZZ_q$, and a complete intersection in $X \times X$.  
We replace the rational parametrization with a local analytic parametrization.
We describe the construction in detail in the application to the explicit 
moduli of abelian varieties of dimensions 1 and 2 described herein, but the 
approach applies in greater generality to any dimension.

\subsection{Complexity hypothesis}
\label{comphyp}

We will denote by $\xf_q$ a finite field of characteristic $p>0$
having $q$ elements. Let $\ZZ_q$ denote the ring of Witt
vectors with values in $\xf_q$. There
exists a canonical lift $\sigma \in \mathrm{Aut}(\ZZ_q)$ of the $p$-th
power Frobenius morphism of $\xf_q$.
If $a$ is an element of $\ZZ_q$ we denote by
$\bar{a}$ its reduction modulo $p$ in $\xf_q$.
We say that we have computed an element $x \in \ZZ_q$ to precision $m$
if we can write down a bit-string representing its class in the
quotient ring $\ZZ_q \slash p^m \ZZ_q$.  In order to assess the
complexity of our algorithms we use the computational model of a
Random Access Machine \cite{MR1251285}. We assume that the
multiplication of two $n$-bit length integers takes $O(n^\mu)$ bit
operations.  One has $\mu=1+\epsilon$ (for $n$ sufficiently large),
$\mu=\log_2 (3)$ and $\mu=2$ using the FFT multiplication algorithm,
the Karatsuba algorithm and a naive multiplication method,
respectively.  Let $x,y \in \ZZ_q \slash p^m \ZZ_q$. For the
following we assume the sparse modulus representation which is explained in
\cite[pp.239]{MR2162716}. Under this assumption one can compute the
product $xy$ to precision $m$ by performing $O(m^\mu \log (q)^\mu)$
bit operations.

\subsection{A lifting algorithm for moduli of elliptic curves}
\label{lift_g1}

We first describe a canonical lifting algorithm for theta null points of 
elliptic curves, hence take an abelian scheme $E$ of relative dimension~$1$ over $\xz_q$. 
Its theta null point $(a_0:a_1:a_2:a_1)$ determines a Legendre model for $E$ 
of the form
$$
y^2 = x(x-1)(x-\lambda), \mbox{ where }
\lambda = \left(\frac{2a_0a_2}{a_0^2+a_2^2}\right)^2.  
$$
In particular we make use of the maps of modular curves 
$$
\cA_1(\Theta_4) \longrightarrow \cA_1(\Theta_4[2]) \longrightarrow X(2),
$$
where the first map is $(a_0:a_1:a_2) \mapsto (a_0:a_2)$ is the restriction 
to the $2$-torsion part of the theta structure, and $X(2)$ is the full 
modular curve of level $2$ with function field generated by $\lambda$.

We recall that the curve $\cA_1(\Theta_4)$ is determined by Riemann's 
equation~\eqref{eqn_riemann_g1} and the correspondence 
equation~\eqref{eqn_corresp_g1} determines a curve in the product 
$\cA_1(\Theta_4) \times \cA_1(\Theta_4)$.  Projecting this correspondence 
curve onto the $2$-torsion part with coordinates $(x_0:x_2)$ and 
$(y_0:y_2)$, gives rise to an affine curve 
\begin{equation}
\label{eqn_corresp_elim_g1}
x^4 - 4x^3 y^3 + 6x^2 y^2 - 4x y + y^4 = 0,
\end{equation}
by setting $x = x_2/x_0$ and $y = y_2/y_0$.  This curve is singular of 
geometric genus $3$, with singularities
$$
\{(0,0), (1,1), (-1,-1), (i,-i)\},
$$ 
where $i^2 = -1$ in $\xz_q$. It is easily verified that all $x$ in 
$\{\infty,0,1,-1,i,i\}$ determine degenerate, singular cubic curves.  
Moreover, the special fiber at $3$ takes the form
$$
(x^3-y)(x-y^3) = 0,
$$ 
whose singularities consist of all points $(x,x^\sigma)$ for $x$ in $\FF_9$.  
Outside of the image of the above degenerate points, the remaining $\FF_9$-rational 
points are supersingular.  

The remaining points correspond to theta null points of ordinary elliptic 
curves, for which it is easily verified that the conditions of 
Lercer-Lubicz~\cite{ll05} for an Artin-Schreier equation are satisfied.  
Hence their Newton algorithm applies to uniquely lift a solution to 
equation~\eqref{eqn_corresp_elim_g1} with the constraint to $y = x^\sigma$.  
From a solution to this system, we set $(a_0:a_2) = (1:x)$ and determine 
$a_1$ by one Newton lifting step.  This gives the following theorem.

\begin{theorem}
\label{existalgo_g1}
There exists a deterministic algorithm which has as input the theta null 
point $(\bar{a}_{i})$ of an elliptic curve $\bar{E}$ over $\xf_q$ and as output the 
theta null point $(a_{i})$ of its canonical lift $E$ to a given precision 
$m \geq 1$, with time complexity 
$$
O(\log(m) d^{\mu}m^\mu)
$$
where $d=\log(q)$.
\end{theorem}

\subsection{A lifting algorithm for split abelian surfaces}
\label{lift_g2_split}

As in Section \ref{equations} we let $\cA_g(\Theta_4[2])$ denote the moduli space of 
$4$-theta null points, projected on the coordinates which are parametrized by
the $2$-torsion subgroup.
We recall that 
$$
a_{00} a_{22} - a_{02} a_{20} = 0,
$$
determines one component in $\cA_2(\Theta_4[2])$ of split abelian 
surfaces.  We refer to Runge~\cite{Runge}) for a complex analytic 
description of this locus as a degenerate Humbert surface.

The remaining components are obtained by the action of a geometric 
automorphism group acting on theta structures and preserving the 
moduli of abelian varieties.  Explicitly this group is generated by 
the projective automorphism group generated by the matrices
$$
\left(\begin{array}{cccc}
1 & 1 & 0 & 0 \\ -1 & 1 & 0 & 0 \\ 0 & 0 & 1 & 1 \\ 0 & 0 & -1 & 1
\end{array}\right), \quad
\left(\begin{array}{cccc}
1 & 0 & 1 & 0 \\ -1 & 0 & 1 & 0 \\ 0 & 1 & 0 & 1 \\ 0 & -1 & 0 & 1
\end{array}\right), \quad
\left(\begin{array}{cccc}
1 & 0 & 0 & 0 \\ 0 & 1 & 0 & 0 \\ 0 & 0 & i & 0 \\ 0 & 0 & 0 & i
\end{array}\right), \quad
\left(\begin{array}{cccc}
1 & 0 & 0 & 0 \\ 0 & i & 0 & 0 \\ 0 & 0 & i & 0 \\ 0 & 0 & 0 & 1
\end{array}\right)\cdot
$$
acting on $\cA_g(\Theta_4[2]) \cong \PP^3$.  
\ignore{
These automorphism determine a group extension $G$ 
$$
1 \rightarrow (\ZZ/2\ZZ) \times (\mu_2)^2 
  \rightarrow G \rightarrow \mathrm{Sp_4(\FF_2)} \rightarrow 1,
$$}
These automorphism determine a transitive action on the $10$ components 
of the Humbert surface.  In particular, given a theta null point of a 
split abelian variety, by means of an automorphism (defined over an 
extension of degree at most $2$), we may assume that it lies on the 
locus $a_{00} a_{22} = a_{02} a_{20}$.

We now recall that the locus $a_{00} a_{22} = a_{02} a_{20}$ is the image 
of $\cA_1(\Theta_4[2]) \times \cA_1(\Theta_4[2])$ in $\cA_2(\Theta_4[2])$ 
by a Segre embedding 
$$
\big((a_0 : a_2),(a_0':a_2')\big) \longmapsto 
(a_{00} : a_{02} : a_{20} : a_{22}) = 
(a_0 a_0' : a_0 a_2' : a_2 a_0' : a_2 a_2').
$$
The canonical lift of this theta null point is obtained by means of 
the canonical lifting to algorithm applied to $(a_{00} : a_{20})
= (a_{02} : a_{22})$ and to $(a_{00} : a_{02}) = (a_{20} : a_{22})$.
This yields the canonical lift of the theta null point with the 
same complexity as for elliptic curves.  We summarize this result 
in the general theorem for abelian surfaces in the next section.

\subsection{A lifting algorithm for moduli of abelian surfaces}
\label{lift_g2}

We use the notation introduced in Section \ref{comphyp}.
For the rest of this section let $A$ be an abelian scheme of relative 
dimension $2$ over $\ZZ_q$ having ordinary reduction.  Suppose $A$ is 
the canonical lift of $A_{\xf_q}$.
Let $\pol$ be an ample symmetric line bundle of degree $1$ on $A$ and
assume we are given a theta structure of type $(\ZZ/4\ZZ)^2$ for
$(A,\pol^4)$.  We denote the theta null point with respect to the
latter theta structure by $(a_{ij})$ where $(i,j) \in (\ZZ/4\ZZ)^2$.
\begin{theorem}
\label{existalgo_g2}
There exists a deterministic algorithm which has as input the theta null 
point $(\bar{a}_{ij})$ of $A_{\xf_q}$ and as output the theta null point 
$(a_{ij})$ of $A$ to a given precision $m \geq 1$, with time complexity 
$$
O(\log(m) d^{\mu}m^\mu)
$$
where $d=\log(q)$.
\end{theorem}

\begin{proof}
The complexity result of Theorem \ref{existalgo_g2} is an analytic version 
of~\cite[Th.2]{ll05}. We explain below how our system of equations can be 
adapted to an analytic context from which the result will follow.
\newline\indent
\ignore{
We embed the open subscheme of theta null points $(a_{ij})$ of abelian surfaces, 
for which $a_{00}$ is invertible, in an affine space, via
\[
(x_{ij}) \mapsto \Big(
\frac{x_{02}}{x_{00}},
\frac{x_{20}}{x_{00}},
\frac{x_{22}}{x_{00}},
\frac{x_{01}}{x_{00}},
\frac{x_{21}}{x_{00}},
\frac{x_{10}}{x_{00}},
\frac{x_{12}}{x_{00}},
\frac{x_{11}}{x_{00}},
\frac{x_{13}}{x_{00}}
\Big)
\]
Let $X$ denote the subset of $\AA^9(\ZZ_q) = \ZZ_q^9$ which is the image of all 
ordinary abelian surfaces under the above defined map.
The Riemann equations (compare Section \ref{equations}) make $X$ an open subset 
of an analytic subspace of $\ZZ_q^9$.
We denote the image of the theta null point of the canonical lift $(a_{ij})$ in 
$X$ by $\alpha$.
The projection $\Psi:X \rightarrow \ZZ_q^3$ to the first three coordinates gives 
a dominant map.  We set $a=\Psi(\alpha)$.} 
By means of a geometric automorphism, we may assume that $a_{00}$ is a unit  and 
embed the corresponding open subscheme of $\cA_2(\Theta_4)$ in $\AA^9$. We identify 
its image in $\cA_2(\Theta_4[2])$ with $\AA^3$.
We denote the open analytic subspace of sections in $\AA^9(\ZZ_q) = \ZZ_q^9$ 
by $X$, and suppose that $\alpha$ is a point of $X$.  This determines a projection 
$\Psi: X \rightarrow \ZZ_q^3$, under which we denote $a = \Psi(\alpha)$.

In the following we let $U \subseteq \ZZ_q^3$ be an analytic neighborhood of $a$, 
and we construct an analytic map $\Phi:U \rightarrow \ZZ_q^3$, such that $\Phi(a) = 0$.  
We first choose pairwise distinct polynomials 
$$
f_1,f_2,f_3 \in \ZZ[\{x_{ij}\},\{y_{ij}\}]
$$ 
from the equations~\eqref{CorrEqns} of Corollary~\ref{genus2}, and let $\Xi$ be the 
function $X \times X \rightarrow \ZZ_q^3$ given by
\begin{eqnarray*}
(x,y)
\mapsto
\big( f_1(x,y),f_2(x,y),f_3(x,y) \big).
\end{eqnarray*}
By the smoothness of the ordinary locus at the prime $3$ of the moduli space 
of abelian surfaces with symmetric theta structure of type $(\ZZ / 4 \ZZ)^2$ 
we conclude that there exists an analytic local inverse $\Pi:U \rightarrow X$ 
of $\Psi$ such that $\Pi(a)=\alpha$ where $U \subseteq \ZZ_q^3$ is a neighborhood 
of $a$ with respect to the $3$-adic topology.  Note that for an arbitrary choice 
of square roots we have
$$
\begin{array}{l@{}l}
a_{01} = \displaystyle\frac{\lambda}{2}\textstyle\!
         \left(\sqrt{b_{00}b_{01}+b_{10}b_{11}}+\sqrt{b_{00}b_{01}-b_{10}b_{11}}\right)\!, \quad &
a_{21} = \displaystyle\frac{\lambda}{2}\textstyle\!
         \left(\sqrt{b_{00}b_{01}+b_{10}b_{11}}-\sqrt{b_{00}b_{01}-b_{10}b_{11}}\right)\!, \\
a_{10} = \displaystyle\frac{\lambda}{2}\textstyle\!
         \left(\sqrt{b_{00}b_{10}+b_{01}b_{11}}+\sqrt{b_{00}b_{10}-b_{01}b_{11}}\right)\!, \quad &
a_{12} = \displaystyle\frac{\lambda}{2}\textstyle\!
         \left(\sqrt{b_{00}b_{10}+b_{01}b_{11}}-\sqrt{b_{00}b_{10}-b_{01}b_{11}}\right)\!, \\ 
a_{11} = \displaystyle\frac{\lambda}{2}\textstyle\! 
         \left(\sqrt{b_{00}b_{11}+b_{01}b_{10}}+\sqrt{b_{00}b_{11}-b_{01}b_{10}}\right)\!, \quad &
a_{13} = \displaystyle\frac{\lambda}{2}\textstyle\! 
         \left(\sqrt{b_{00}b_{11}+b_{01}b_{10}}-\sqrt{b_{00}b_{11}-b_{01}b_{10}}\right)\!,
\end{array}
$$
where
$$
\begin{array}{ll}
b_{00} = 1, \quad &
b_{01} = \sqrt{\lambda^{-1}(a_{00}a_{02}+a_{20}a_{22})}, \\
b_{20} = \sqrt{\lambda^{-1}(a_{00}a_{20}+a_{02}a_{22})}, \quad &
b_{22} = \sqrt{\lambda^{-1}(a_{00}a_{22}+a_{02}a_{20})},
\end{array}
$$
\ignore{DRK: I prefer the above notation b_ij which suggests their role as (Z/2Z)^2-theta null points.
$$
\begin{array}{ll}
a_{01} = \lambda(\sqrt{r_{1}r_{2}+r_{3}r_{4}}+\sqrt{r_{1}r_{2}-r_{3}r_{4}})/2, \quad &
a_{21} = \lambda(\sqrt{r_{1}r_{2}+r_{3}r_{4}}-\sqrt{r_{1}r_{2}-r_{3}r_{4}})/2 \\
a_{10} = \lambda(\sqrt{r_{1}r_{3}+r_{2}r_{4}}+\sqrt{r_{1}r_{3}-r_{2}r_{4}})/2, \quad &
a_{12} = \lambda(\sqrt{r_{1}r_{3}+r_{2}r_{4}}-\sqrt{r_{1}r_{3}-r_{2}r_{4}})/2 \\
a_{11} = \lambda(\sqrt{r_{1}r_{4}+r_{2}r_{3}}+\sqrt{r_{1}r_{4}-r_{2}r_{3}})/2, \quad &
a_{13} = \lambda(\sqrt{r_{1}r_{4}+r_{2}r_{3}}-\sqrt{r_{1}r_{4}-r_{2}r_{3}})/2,
\end{array}
$$
where
$$
\begin{array}{ll}
r_{1} = 1, \quad &
r_{2} = \sqrt{\lambda^{-1}(a_{00}a_{02}+a_{20}a_{22})}, \\
r_{3} = \sqrt{\lambda^{-1}(a_{00}a_{20}+a_{02}a_{22})}, \quad &
r_{4} = \sqrt{\lambda^{-1}(a_{00}a_{22}+a_{02}a_{20})},
\end{array}
$$}
and where $\lambda = (a_{00}^2+a_{02}^2+a_{20}^2+a_{22}^2)/2$.  
The above formulas can be deduced from Mumford's $2$-multiplication
formula \cite[$\S$3]{mu66}.
We note that the zero set of $b_{ij}$ and of
$$
b_{00}b_{01} \pm b_{10}b_{11}, \quad
b_{00}b_{10} \pm b_{01}b_{11}, \quad 
b_{00}b_{11} \pm b_{01}b_{10},
$$
lie over the components of the moduli of split abelian varieties.  Applying 
the algorithm of the previous section to such points, we may thus assume that 
the map is unramified at $(a_{ij})$.
\newline\indent
Consider the subset $Y \subseteq X \times X$ which is defined by the equations
of Corollary~\ref{genus2}.
Let $p_i:Y \subseteq X \times X \rightarrow X$ ($i=1,2$) be the map
induced by the projection on the $i$th factor.
The map $p_i$ forms an \'{e}tale covering.
We can choose a local analytic section $i_1:V \rightarrow Y$ of the projection $p_1$
in a neighbourhood $V$ of $\alpha$ such that $i_1(\alpha)=(\alpha,\alpha^{\sigma})$.
Let $\Sigma = p_2 \circ i_1$.
Note that $\Sigma(\alpha)=\alpha^{\sigma}$.
By Serre-Tate theory the morphism $\Sigma$ is analytic
in a neighborhood of $\alpha$.
More precisely, this is a consequence of
the representability of the local deformation space of an ordinary
abelian variety over $\xf_q$ by a formal torus and the fact that the
unique lift of the relative $3$-Frobenius equals up to isomorphism
the $3$rd powering map on this torus (see \cite{ka81} and \cite{mo95}).
\newline\indent 
We define $\Phi$ to be the composition
$$
\SelectTips{cm}{10}
\xymatrix@C=0.75cm{
U \ar[r]^-{\Delta} &
U \times U \ar[r]^-{\Pi^2} &
X \times X \ar[r]^-{\mathrm{Id} \times \Sigma} &
X \times X \ar[r]^-{\Xi} & \ZZ_q^3
}
$$
where $\Delta$ is the diagonal map and $\Pi^2= \Pi \times \Pi$.
The equality $\Phi(a) = 0$ holds by Corollary~\ref{genus2}.
By the above discussion, $\Phi$ is analytic and is defined on the
open disc $U$ with center $a$ and radius $1$. The fact that the radius
equals $1$ can be deduced from the interpretation of the points
in the image of $\Psi$ as moduli points of abelian varieties with $2$-theta structure.
\newline\indent
In the following we verify that the assumptions of
\cite[Th.2]{ll05} are satisfied.
For an analytic function $F$ we
denote its first order derivative by $D_{F}$. We have
$\Phi(x) \equiv 0 \bmod 3$ for all $x \in U$, because all points of
$U$ reduce to the same canonical theta null point satisfying the
equations (\ref{RiemEqns}) and (\ref{CorrEqns}), and hence $D_{\Phi}(a)
\equiv 0 \bmod 3$. We write $D_{\Xi_X}$ and $D_{\Xi_Y}$ for the
submatrices of $D_{\Xi}$ being the derivative of $\Xi$ with respect to
the first and second factor of the product $X \times X$.  By the chain
rule we conclude that
\begin{eqnarray}
\label{congu}
0 \equiv D_{\Phi}(a) = D_{\Xi_X} \big( \Pi(a),\Pi(a)^{\sigma}
\big) D_{\Pi}(a)
+ D_{\Xi_Y} \big( \Pi(a),\Pi(a)^{\sigma} \big)
D_{\Sigma} \big( \Pi(a) \big) D_{\Pi}(a) \bmod 3.
\end{eqnarray}
By general theory \cite{ka81}
the Frobenius lift acts on the Serre-Tate formal torus
as the $3$-rd powering map and hence $D_{\Sigma} \big( \Pi(a) \big) \equiv 0 \bmod 3$.
We conclude by equation (\ref{congu}) that
\begin{eqnarray}
\label{condition1}
D_{\Xi_X} \big( \Pi(a),\Pi(a)^{\sigma}
\big) D_{\Pi}(a) \equiv 0 \bmod 3.
\end{eqnarray}
\indent Next we prove by contradiction that for a suitable choice of the triple
$(f_1,f_2,f_3)$ (notation as above) we have
\begin{eqnarray}
\label{condition2}
\mathrm{det} \Big( D_{\Xi_Y} \big( \Pi(a),\Pi(a)^{\sigma}
\big) D_{\Pi}(a) \Big) \not\equiv 0 \bmod 3.
\end{eqnarray}
We remark that in the lifting algorithm the choice of the triple
$(f_1,f_2,f_3)$ has to be done depending on the initial data.
Suppose condition (\ref{condition2}) is not satisfied for any triple $(f_1,f_2,f_3)$.
Then by the Jacobi criterion we conclude
that the moduli space of pairs of ordinary abelian surfaces with symmetric
$4$-theta structure and compatible $(3,3)$-isogeny is not
smooth at the point $\big( \Pi(a),\Pi(a)^{\sigma} \big)$.
This contradicts the fact that the latter space forms an
\'{e}tale covering of the smooth space $X$.
Clearly the equations (\ref{condition1}) and (\ref{condition2}) imply the
assumptions of \cite[Th.2]{ll05}.
By the algorithm suggested there and the above discussion
we can compute $x \in U$
such that $\Phi(x) \equiv 0 \bmod 3^m$ for given precision $m$
with complexity as stated in the theorem.
\newline\indent
In the following we explain why the output of the latter algorithm
is indeed the theta null point of the
canonical lift to given precision.
We claim that for $x \in U$ one has an equivalence
\begin{eqnarray}
\label{equivalence}
x \equiv a \bmod 3^m \quad \Leftrightarrow \quad \Phi(x) \equiv 0
\bmod 3^m.
\end{eqnarray}
It suffices to prove that $\Phi(x) \equiv 0
\bmod 3^m$ implies $x \equiv a \bmod 3^m$ since the
converse is obvious.
The proof is done by induction on $m$.
Assume that equivalence (\ref{equivalence}) holds for
some $m \geq 1$.
Assume that $\Phi(x) \equiv 0 \bmod 3^{2m}$.
By the induction hypothesis we have $\delta = 3^{-m}(x-a) \in \ZZ_q^3$.
Then by Taylor expansion of the analytic function $\Phi$ at $a$ we get
\begin{eqnarray}
\label{reform}
0 \equiv \Phi(x) =\Phi(a + 3^m \delta) \equiv \Phi(a) + D_{\Phi}(a)
3^m \delta + \ldots \equiv D_{\Phi}(a) 3^m \delta \bmod 3^{2m}.
\end{eqnarray}
By equation (\ref{reform}) we conclude that
\begin{eqnarray}
\label{reform2}
0 \equiv D_{\Phi}(a) \delta \bmod 3^m.
\end{eqnarray}
We set
\begin{eqnarray*}
D_X & = & D_{\Xi_X} \big( \Pi(a),\Pi(a^{\sigma}) \big) D_{\Pi}(a), \\
D_Y & = & D_{\Xi_Y} \big( \Pi(a),\Pi(a^{\sigma}) \big) D_{\Pi}(a).
\end{eqnarray*}
By (\ref{congu}) the equation (\ref{reform2}) is equivalent to
\begin{eqnarray}
\label{fixpoint}
\delta \equiv D(\delta) \bmod 3^m
\end{eqnarray}
where $D$ is the linear operator given by
\[
y \mapsto -\left( D_Y^{-1} D_X y \right)^{\sigma^{-1}}
\]
Here we have used that the map $\Sigma$ already exists as an
endomorphism of $U$ which commutes with the application $\Pi$.
Note that by condition
(\ref{condition2}) the matrix $D_Y$ is invertible.  By condition
(\ref{condition1}) the entries of the matrix $D_Y^{-1} D_X$ are all
divisible by $3$.
As a consequence we conclude from equation (\ref{fixpoint})
that $\delta \equiv 0 \bmod 3^m$.
This proves our claim.
\newline\indent
In the following we will show how to compute the matrices $D_X$ and $D_Y$,
since they are needed for the algorithm of Lercier and Lubicz~\cite{ll05}.
By the above discussion, we can compute compatible branches of the local 
inverse $\Pi$ at $a$ and $a^\sigma$ such that $\Pi(a^{\sigma}) = 
\Pi(a)^{\sigma}$.
From this it is straightforward to compute
$D_{\Xi_X} \big( \Pi(a),\Pi(a^{\sigma}) \big)$ and $D_{\Xi_Y} \big(
\Pi(a),\Pi(a^{\sigma}) \big)$.
Next we explain how to compute $D_{\Pi}(a)$.
Let $\Lambda:\ZZ_q^9 \rightarrow \ZZ_q^{20}$ be defined by
$$
x = (x_{ij}) \longmapsto (\Lambda_1(x),\dots,\Lambda_{20}(x))
$$
where $\Lambda_i$ are the Riemann relations~\eqref{RiemEqns}, 
so that $\Lambda(\Pi(a)) = 0$. By the chain rule we conclude that
\begin{eqnarray}
\label{implicit}
D_{\Lambda} \big( \Pi(a) \big) D_{\Pi}(a) = 0.
\end{eqnarray}
Let $\pi:\ZZ_q^3 \rightarrow \ZZ_q^6$ be the morphism such
that $\Pi(a)= \big( a,\pi(a) \big)$.
Then $D_{\Pi}$ is the vertical join of the unit matrix of rank $3$
and $D_{\pi}$ where $D_{\pi}$ denotes the derivative of $\pi$.
We write
\[
D_{\Lambda}= \left( D_{\Lambda}^{(1)}  D_{\Lambda}^{(2)} \right)
\]
where $D_{\Lambda}^{(1)}$ and $D_{\Lambda}^{(2)}$ have $3$ and $6$ columns,
respectively.
Note that by the smoothness of the space $X$ the rank of
$D_{\Lambda}^{(2)}$ at $\Pi(a)$ equals $6$. There exists a
matrix $T \in \mathrm{GL}_{20}(\ZZ_q)$ such that the matrix
\[
E=T \cdot D_{\Lambda}^{(2)} \big( \Pi(a) \big)
\]
is in echelon form.
It follows from equation (\ref{implicit}) and the above discussion that
\begin{eqnarray*}
E \cdot D_{\pi}(a) = - T \cdot D_{\Lambda}^{(1)} \big( \Pi(a) \big).
\end{eqnarray*}
From this it is straightforward to compute $D_{\pi}(a)$ inverting the unique 
invertible $(6,6)$-square submatrix of $E$.  We remark that the above computation 
can be done modulo any given precision. 
This completes the proof of Theorem \ref{existalgo_g2}.
\end{proof}

We conclude this section by a practical remark.  Our implementation uses a 
multivariate version of the algorithm of R. Harley (compare \cite[$\S$3.10]{Vercauteren2003})
for solving generalized Artin-Schreier equations instead of the the method 
suggested in \cite{ll05}.

\subsection{LLL reconstruction}
\label{LLL}

From the theory of complex multiplication we know that the invariants 
of canonical lifts are algebraic over $\QQ$.  We briefly recall the 
method of Gaudry et al.~\cite{Gaudryandall} for LLL reconstruction of 
algebraic relations over $\ZZ$.  Let $\gamma$ be a $p$-adic integer in an 
extension of degree $r$ over $\ZZ_p$, and let $m$ be the precision to 
which it is determined.  We assume that the degree $n$ of its minimal 
polynomial over $\QQ$ is known, i.e.~that there exists $f(x) \in \ZZ[x]$, 
with
$$   
f(\gamma) = a_n\gamma^n + \hdots + a_0 = 0,
$$   
where the $a_i \in \ZZ$ are unknown. We determine a basis of the left 
kernel in $\ZZ^{n+r+1}$ of the vertical join of the matrix 
$$   
\left[
\begin{array}{cccc}    
1 & 0 & \cdots & 0\\   
\gamma_{1, 0} & \gamma_{1, 1} & \cdots &\gamma_{1, (r - 1)}\\    
\vdots& &  & \vdots\\  
\gamma_{n, 0} & \gamma_{n, 1} & \cdots &\gamma_{n, (r - 1)}\\    
\end{array}
\right]    
$$   
with $p^m$ times the $r\times r$ identity matrix, where 
$\gamma_{i, j}$ are defined by    
$$    
\gamma^i = \gamma_{i,0} + \gamma_{i,1}w_1 + \hdots + \gamma_{i, (r - 1)}w_{r-1},
$$   
in terms of a $\ZZ_p$ basis $\{1,w_1,\dots,w_{r-1}\}$ for $\ZZ_q$.
The minimal polynomial $f(x)$ is determined by LLL as a short vector
$(a_0, \hdots, a_n, \varepsilon_1, \hdots, \varepsilon_r)$ in the 
kernel. 
\newline\indent
The complexity of the LLL step depends on the values $r$, $n$, and $m$. 
The values of $r$ and $n$ can be recovered by a curve selection and analysis of 
the Galois theory of the class fields.  The required precision $m$, 
determined by the size of the output, is less well-understood, and we 
express the complexity in terms of these three parameters.
Using the $L^2$ variant of LLL by Nguy$\tilde{\mbox{\^ e}}$n and 
Stehl{\'e}~\cite{NgSt05}, the complexity estimate of~\cite{Gaudryandall} 
gives $O((n+r)^5(n+r+m)m)$ in general, and in our case the special 
structure of the lattice gives a complexity of $O((n+r)^4(n+r+m)m)$.

\section{Moduli equations and parametrizations}
\label{riemequations}

In this section we give the equations which form a higher dimensional 
analogue of Riemann's quartic theta relation.  Then we state the classical 
Thomae formulas in an algebraic context, relating the invariants of 
genus~$2$ curves to theta null points.  Finally we apply the algorithm 
of Section~\ref{algo} to the construction of CM invariants of abelian 
surfaces and genus 2 curves.

\subsection{The Thomae formulas for genus $2$}
\label{thomae}

Let $R$ be an unramified local ring of odd residue characteristic, and 
$H$ a hyperelliptic curve over $R$ given by an affine equation
\[
y^2=\prod_{i=1}^5 (x-e_i),
\]
where the $e_i \in R$ are pairwise distinct in the residue field.
Let $(J,\varphi)$ denote the Jacobian of $H$ where $\varphi$ denotes the 
canonical polarization.
There exists a finite unramified extension $S$ of $R$ and an ample symmetric 
line bundle $\pol$ of degree $1$ on $J_S$ which induces the polarization 
$\varphi_S$.
We assume that $S$ is chosen such that there exists an $S$-rational symmetric 
theta structure $\Theta$ of type $(\xz/4 \xz)^2$ for the pair $(J_S,\pol)$.
Let $(a_{ij})$, where $(i,j) \in (\xz/4 \xz)^2$, denote the theta null point 
with respect to the latter theta structure.
\begin{theorem}[Thomae formulas]
\label{tthomae}
With the notation as above, one has
$$
\begin{array}{ll}
a_{00} = 1 & 
\displaystyle
a_{02} = \sqrt[4]{
\frac{(e_1-e_4)(e_2-e_5)(e_3-e_4)}{(e_1-e_5)(e_2-e_4)(e_3-e_5)}} \\
\displaystyle
a_{20} = \sqrt[4]{
\frac{(e_1-e_2)(e_1-e_4)}{(e_1-e_3)(e_1-e_5)}} &
\displaystyle
a_{22} = \sqrt[4]{
\frac{(e_1-e_2)(e_2-e_5)(e_3-e_4)}{(e_1-e_3)(e_2-e_4)(e_3-e_5)}}
\end{array}
$$
such that $a_{02}^2 = (e_1 - e_3)/(e_1 - e_2)a_{20}^2a_{22}^2$.
\end{theorem}

\noindent
For a complex analytic proof of the Thomae formulas see~\cite[p.120]{mu84}.
\newline\indent
Conversely, let $A$ be an abelian surface over $S$ with ample symmetric 
line bundle $\pol$ of degree $1$ on $A$. 
Assume we are given a symmetric theta structure of type $(\xz/4\xz)^2$ 
for the pair $(A,\pol^4)$, and let $(a_{ij})$ denote the theta null point 
with respect to the latter theta structure.  We associate a curve to the 
theta null point in the following way. Let $\mu$ be a solution of the 
equation (possibly over an unramified extension)
$$
\mu^2 - 
\frac{(a_{00}^4-a_{02}^4+a_{20}^4-a_{22}^4)}
{(a_{00}^2a_{20}^2-a_{02}^2a_{22}^2)}\mu + 1 = 0,
$$
and set 
$$
\lambda_1 = \Big(\frac{a_{00}a_{02}}{a_{22}a_{20}}\Big)^2\!\!, \quad
\lambda_2 = \Big(\frac{a_{02}}{a_{22}}\Big)^2\!\!\mu, \quad
\lambda_3 = \Big(\frac{a_{00}}{a_{20}}\Big)^2\!\!\mu\cdot
$$
\begin{corollary}
The curve
\[
y^2 = x(x-1)(x-\lambda_1)(x-\lambda_2)(x-\lambda_3)
\]
has as Jacobian the abelian surface $A$.
\end{corollary}

\begin{proof}
Inverting Theorem~\ref{tthomae}, one verifies that the roots $\mu$ give 
rise to values $(\lambda_1,\lambda_2,\lambda_3)$ in the set
$$
\left\{
\left(
    \frac{e_1 - e_3}{e_1 - e_2},
    \frac{e_1 - e_4}{e_1 - e_2},
    \frac{e_1 - e_5}{e_1 - e_2}\right), 
\left(
    \frac{e_1 - e_3}{e_1 - e_2},
    \frac{e_1 - e_3}{e_1 - e_5},
    \frac{e_1 - e_3}{e_1 - e_4}\right)
\right\},
$$
which determine curves isomorphic to the curve with affine equation 
$
y^2 = \prod_{i=1}^{5} (x-e_i).
$
\end{proof}

Unfortunately, the Rosenhain invariants $(\lambda_1,\lambda_2,\lambda_3)$ 
of the above curve are not determined by the $2$-torsion part 
$(a_{00}:a_{02}:a_{20}:a_{22})$ of the theta null point.  Instead we must 
pass to a $(2,2)$-isogenous curve to determine a genus 2 curve parametrized 
by this theta null point.

\begin{theorem}
The curve
$$
y^2 = x(x-1)(x-\mu_1)(x-\mu_2)(x-\mu_3),
$$
where
$$
\begin{array}{l}
\mu_1 = \displaystyle 
\frac{(a_{00}^2+a_{02}^2+a_{20}^2+a_{22}^2)(a_{00}a_{02}+a_{20}a_{22})}%
     {2(a_{00}a_{20}+a_{02}a_{22})(a_{00}a_{22}+a_{02}a_{20})} \\ \\

\mu_2 = \displaystyle  
\frac{(a_{00}^2-a_{02}^2+a_{20}^2-a_{22}^2)(a_{00}a_{02}+a_{20}a_{22})}%
     {2(a_{00}a_{22}+a_{02}a_{20})(a_{00}a_{20}-a_{02}a_{22})}\\ \\
\mu_3 = \displaystyle
\frac{(a_{00}^2+a_{02}^2+a_{20}^2+a_{22}^2)(a_{00}^2-a_{02}^2+a_{20}^2-a_{22}^2)}%
     {(a_{00}a_{02}+a_{20}a_{22})(a_{00}a_{02}-a_{20}a_{22})}
\end{array}
$$
has Jacobian $(2,2)$-isogenous of the abelian surface $A$.
\end{theorem}

\begin{proof}
The Richelot isogeny determined by the polynomials 
$$
G_1 =  x, \quad 
G_2 = (x-1)(x-\lambda_1), \quad 
G_3 = (x-\lambda_2)(x-\lambda_3),
$$
determines a curve isomorphic to the above curve.
\end{proof}

Thus we obtain a rational map from the space $\cA_2(\Theta_4[2])$, 
determined by the $2$-torsion part $(a_{00}:a_{02}:a_{20}:a_{22})$ 
of a theta null point, to the moduli space $\cM_2(2)$ of genus $2$ 
curves with level-$2$ structure, determined by the Rosenhain invariants 
$(\mu_1,\mu_2,\mu_3)$.  The latter point specifies an ordered 
six-tuple of Weierstrass points over $(\infty,0,1,\mu_1,\mu_2,\mu_3)$.
We note that this map is defined on the open subspace outside of 
the components defining split abelian surfaces.

\subsection{Examples of canonical lifts}
\label{example}

In this section we give examples of canonical lifts of $3$-adic theta null 
points. The examples were computed using implementations of our algorithms 
in the computer algebra system Magma~\cite{Magma}.  Generic algorithms and 
databases of CM invariants for genus 2 curves can be found from the authors'
web pages (see~\cite{KohelCode}).

\medskip

\noindent{\bf Example 1.}
Consider the genus 2 hyperelliptic curve $\bar{H}$ over $\xf_3$ defined by 
the equation
\[
y^2=x^5+x^3+x+1.
\]
Let $\bar{J}$ denote the Jacobian of $\bar{H}$.
The abelian surface $\bar{J}$ is ordinary.
Over an extension of degree $40$
there exists a theta structure of type $(\xz / 4 \xz)^2$
for $(\bar{J},\pol^4)$ where $\pol$ is the line bundle corresponding
to the canonical polarization.
Let $(\bar{a}_{ij})$ denote the theta null point of $(\bar{J},\pol^4)$
with respect to the latter theta structure.
We can assume that $\bar{a}_{00}=1$.
Note that the coordinates $\bar{a}_{02}$, $\bar{a}_{20}$ and $\bar{a}_{22}$
are defined over an extension of degree $10$.
We set $\xf_{3^{10}}=\xf_3[z]$ where 
$z^{10} + 2z^6 + 2z^5 + 2z^4 + z + 2 = 0$.
We choose
\[
\bar{a}_{02}=z^{9089}, \quad \bar{a}_{20}=z^{18300} \quad \mbox{and}
\quad \bar{a}_{22}=z^{8601}.
\]
By the algorithm described in Section \ref{algo} we lift the
triple $(\bar{a}_{02},\bar{a}_{20},\bar{a}_{22})$
to the unramified extension of $\xz_3$ of degree $10$.
We denote the lifted coordinates by
$a_{02}$, $a_{20}$ and
$a_{22}$.
Let $P_{ij}$ be the minimal polynomial of $a_{ij}$
over $\xq$.
A search for algebraic relations using the LLL-algorithm yields 
\begin{eqnarray*}
P_{02} & = & x^{80} - 69x^{76} + 4911x^{72} + 20749x^{68} +
299094x^{64}  - 202217x^{60} \\
& & + 1093161x^{56} - 7393871x^{52} + 11951456x^{48} + 7541235x^{44}
\\
& & - 26349059x^{40} + 7541235x^{36} + 11951456x^{32} - 7393871x^{28}
\\
& & + 1093161x^{24} - 202217x^{20} + 299094x^{16} + 20749x^{12} +
4911x^8 \\
& & - 69x^4 + 1, \\
P_{20} & = & x^{20} - 5x^{19} + 23x^{18} - 53x^{17}
+ 112x^{16} - 203x^{15} + 279x^{14} - 345x^{13} \\
& & + 360x^{12} - 333x^{11} + 329x^{10} - 333x^9
+ 360x^8 - 345x^7 + 279x^6 \\
& & - 203x^5 + 112x^4 - 53x^3 + 23x^2 - 5x + 1, \\
P_{22} & = & x^{80} + 5x^{76} + 184x^{72} + 2254x^{68}
+ 4470x^{64} + 160109x^{60} + 768428x^{56} \\
& & + 421488x^{52} + 36971535x^{48} - 75225290x^{44} + 44767882x^{40}
\\
& & - 43287046x^{36} + 86078086x^{32} - 75568556x^{28} +
31873762x^{24} \\
&& - 7293064x^{20} + 989181x^{16} - 32859x^{12} + 4318x^8 + 44x^4 + 1.
\end{eqnarray*}
We conclude that the field $k_0$ generated by the coordinates $a_{02}$,
$a_{20}$ and $a_{22}$
is a Galois extension of $\xq$ having degree $160$. Note that $k_0$ contains $\xq(i)$.
\newline\indent
The characteristic polynomial of the absolute Frobenius endomorphism of
$\bar{J}$ equals
\[
x^4 + 3 x^3 + 5  x^2 + x + 9.
\]
Let $K = \mathrm{End}_{\xf_3}(\bar{J}) \otimes \xq$.
The field $K$ is a normal CM field of dimension $4$ whose Galois group equals
$\xz / 4 \xz$. The class number of $K$ equals $1$.
The maximal totally real subfield of $K$ is given by $\xq(\sqrt{13})$.
Note that $K$ equals its own reflex field $K^*$.
The compositum $k_0K^*$ forms an abelian extension of $K^*$ having
conductor $8$ and Galois group $(\xz / 2 \xz)^2 \times \xz/ 10 \xz$.
Note that the polynomial $P_{20}$ generates the ray class
field of $K^*$ modulo $2$.
\newline\indent
We remark that the curve $H$ with defining equation
$$
y^2 = 52x^5-156x^4+208x^3-156x^2+64x-11
$$
is a canonical lift of $\bar{H}$ in the sense that $H$ reduces to the 
curve $\bar{H}$ and the Jacobian of $H$ is isomorphic to the canonical 
lift of $\bar{J}$.  For a list of curves of genus $2$ over $\xq$ 
having complex multiplication we refer to \cite{vw99}.

\medskip

\noindent{\bf Example 2.}
Let $\bar{H}$ be the hyperelliptic curve over $\FF_{3^6} = \FF_{3}[z]$ where 
$z^3 - z + 1 = 0$, defined by the affine equation
$$
y^2 = x(x-1)(x-z)(x-z^8)(x-z^2).
$$
We may associate a theta null point to $\bar{H}$ over an extension and apply 
our algorithm to determine the canonical lifted Rosenhain invariants from the 
lift of the theta null point.  By LLL reconstruction, the Igusa invariants 
$$
j_1 = \frac{J_2^5}{J_{10}},\quad
j_2 = \frac{J_2^3J_4}{J_{10}},\quad
j_4 = \frac{J_2 J_8}{J_{10}},
$$
of the canonically lifted curve $H$ satisfy the minimal polynomials 
{\scriptsize 
$$
\begin{array}{l}
1167579244112528766379604000052855618647029683j_1^6 \\ \quad -\, 
15257677849803613955571236222133142793627666039890131548110848j_1^5 \\ \quad +\, 
1196131879277094213213237826625656616667290986216439120696238769598103552j_1^4 \\ \quad -\, 
1502690183964538566290599551441994054504503089078463931648679137089316924162048j_1^3 \\ \quad +\, 
9494960051498045134856366244512386171442968847268749046183153757319495998347673600000j_1^2 \\ \quad -\, 
9489242494532768198621993753759532669268063460725268563272920396343489385558179840000000000j_1 \\ \quad +\, 
3154745183558232433309182902721654489400212652045192101874580090073682713333727232000000000000000 \\
31524639591038276692249308001427101703469801441j_2^6 \\ \quad -\, 
16745634807723620828207592940844036495138204085628428409110528j_2^5 \\ \quad -\, 
12265164179615739710029144012197055859859725320474999182497036825001984j_2^4 \\ \quad +\, 
352141775319032803460285640460530428476805227032807841788375367068285927424j_2^3 \\ \quad -\, 
115886117015701373170818041387627276397709556079989081954770457714548434534400000j_2^2 \\ \quad +\, 
6241088101000204747012315559761320786612924621590641411279130896395801722880000000000j_2 \\ \quad -\, 
119116948667007461483450210289814155018636097544277843858343319784591982592000000000000000 \\
22981462261866903708649745533040357141829485250489j_4^6 \\ \quad -\, 
38333133385822330975872342595626396239705000243787196311246336j_4^5 \\ \quad -\, 
13445890564402694049486311582599736771794395285600128293985309687808j_4^4 \\ \quad -\, 
25587083283087299157726904789352095023627415391850896175427316095123456j_4^3 \\ \quad -\, 
20922653078662308982945894934868322119306736601817862795598824527101952000j_4^2 \\ \quad -\, 
6125981423009705673176896782997851830442900916324351082547267950870528000000j_4 \\ \quad -\, 
1226005575547426252457067048464156648937773482166774996185845610840064000000000 \\
\end{array}
$$}
We note that neither of these Jacobians has good ordinary reduction at $2$, thus 
extend the realm of applicability of the $2$-adic CM method~\cite{Gaudryandall}.

\section{Conclusion and perspectives}

This work generalizes prior higher dimensional $2$-adic canonical lifting 
algorithms to a $3$-adic setting.  Firstly, in Theorem~\ref{kop}, we introduce 
the moduli equations which provide the higher dimensional analogues of the 
modular curve $X_0(3)$.  Secondly, we describe a general multivariate Hensel 
lifting algorithm in an analytic framework (removing the need for a rational 
parametrization of a variety).  
As an application our work gives an explicit CM construction for moduli of 
genus 2 curves (and their Jacobian surfaces), yielding a $3$-adic alternative 
to the $2$-adic construction of Gaudry et al.~\cite{Gaudryandall}, and  
extending the domain of applicability to additional quartic CM fields.  
We expect that our method extends to primes $p>3$, for which the primary 
ingredient will be an analogue of our Theorem~\ref{kop}.  
With an increasing complexity for the resulting schemes, as both $p$ and the 
dimension grow, we expect our approach through analytic parametrizations will become 
essential.  

\bigskip

\noindent{\bf Acknowledgments}.
We are grateful to Y. Kopeliovich for explaining his results
and for giving some very valuable references to the
literature.

\bibliographystyle{plain}

\end{document}